\documentclass[a4paper, 11pt]{article}
\usepackage{latexsym,amsfonts, amsmath, amssymb, amsthm}
\usepackage[utf8]{inputenc}
\usepackage[final=true,protrusion=true,expansion=true]{microtype}
\usepackage[noadjust]{cite}
\usepackage{enumitem}
\usepackage{color}
\usepackage{mathtools}
\usepackage{chngcntr}
\usepackage{apptools}
\AtAppendix{\counterwithin{theorem}{section}}

\usepackage{booktabs}
\usepackage{multirow}

\usepackage{subcaption}
\usepackage[font=scriptsize,labelfont=bf]{caption}
\usepackage[affil-it]{authblk}

\usepackage{geometry}
\usepackage{hyperref}
\hypersetup{
	final,
    colorlinks=true,
    linkcolor=blue,
    filecolor=magenta,
    urlcolor=cyan,
}

\theoremstyle{plain}
\newtheorem{theorem}{Theorem}

\newtheorem{lemma}[theorem]{Lemma}

\theoremstyle{plain}
\newtheorem{definition}[theorem]{Definition}

\newtheorem{remark}[theorem]{Remark}

\providecommand{\Supp}{\operatorname{supp}}                            
\providecommand{\supp}{\Supp}

\providecommand{\argmin}{\operatorname*{\arg\min}}  
\providecommand{\argmax}{\operatorname*{argmax}}  
\providecommand{\Id}{\Op{Id}}                     
\providecommand{\diag}{\operatorname{diag}}

\providecommand{\limsup}{\operatorname*{limsup}}

\providecommand{\CC}{{\cal C}}

\providecommand{\CE}{{\cal E}}

\providecommand{\CM}{{\cal M}}
\providecommand{\CN}{{\cal N}}

\providecommand{\CP}{{\cal P}}

\providecommand{\CT}{{\cal T}}

\providecommand{\CX}{{\cal X}}
\providecommand{\CY}{{\cal Y}}

\providecommand{\bbE}{\mathbb{E}}

\providecommand{\bbN}{\mathbb{N}}

\providecommand{\bbR}{\mathbb{R}}

\providecommand*{\abs}[1]{\left|{#1}\right|} 
\providecommand*{\absnormal}[1]{|{#1}|} 
\providecommand*{\N}[1]{\left\|{#1}\right\|} 
\providecommand*{\Nnormal}[1]{\|{#1}\|} 
\providecommand*{\Nbig}[1]{\big\|{#1}\big\|} 

\newcommand*{\SN}[1]{\left|{#1}\right|}      

\providecommand*{\abs}[1]{\left|{#1}\right|} 

\newcommand*{\Op}[1]{\mathsf{#1}} 

\usepackage{bbm}
\usepackage{arydshln}
\usepackage{subcaption}
\usepackage{xcolor}
\usepackage{ifdraft}
\usepackage{subeqnarray}
\usepackage{cases}
\usepackage{bm}

\usepackage{algorithm}
\usepackage{algpseudocode}

\newif\ifrevised
\newcommand{\revised}[1]{%
	\ifrevised
		\color{purple} #1 \color{black} %
	\else
		#1%
	\fi}
\revisedfalse

\newif\ifrevisedTwo

\revisedTwotrue

\newcommand{\overbar}[1]{\mkern 1.5mu\overline{\mkern-1.5mu#1\mkern-1.5mu}\mkern 1.5mu}
\newcommand{\overbarwithsubscript}[1]{\makebox[0pt]{$\phantom{#1}\mkern 1.5mu\overline{\mkern-1.5mu\phantom{#1}\mkern-1.5mu}\mkern 1.5mu$}#1}
\renewcommand{\underbar}[1]{\mkern 1.5mu\underline{\mkern-1.5mu#1\mkern-1.5mu}\mkern 1.5mu}

\newcommand{\overbarscript}[1]{\mkern 1.5mu\overline{\mkern-1.5mu#1\mkern-1.5mu}\mkern 0mu}

\newcommand{\divergence}{\textrm{div}}

\DeclareMathOperator*{\Law}{Law}

\newcommand{\saddlepoint}{x^*,y^*}

\newcommand{\VarX}{\mathrm{Var}^X}
\newcommand{\VarY}{\mathrm{Var}^Y}

\newcommand{\OX}{\overbarwithsubscript{X}}
\newcommand{\OY}{\overbarwithsubscript{Y}}

\newcommand{\omegaa}[0]{\omega_{\alpha}}
\newcommand{\omegab}[0]{\omega_{-\beta}}

\newcommand{\conspointx}[1]{x_{\alpha}^Y({#1})}

\newcommand{\conspointy}[1]{y_{\beta}^X({#1})}

\newcommand{\empmeasure}[1]{\widehat\rho_{#1}^N}
\newcommand{\empmeasureX}[1]{\widehat\rho_{X,#1}^{\revised{N_1}}}
\newcommand{\empmeasureY}[1]{\widehat\rho_{Y,#1}^{\revised{N_2}}}

\newcommand*\samethanks[1][\value{footnote}]{\footnotemark[#1]}

\makeatother

\title{\usefont{OT1}{bch}{b}{n}
	\huge Consensus-Based Optimization\\for Saddle Point Problems \\
}

\date{}

\author[1]{Hui Huang\thanks{Email: \texttt{hui.huang@uni-graz.at}, \texttt{jinniao.qiu@ucalgary.ca}, \texttt{konstantin.riedl@ma.tum.de}}}
\affil[1]{University of Graz, Institute of Mathematics and Scientific Computing, Graz, Austria}
\author[2]{Jinniao Qiu\samethanks[1]}
\affil[2]{University of Calgary, Department of Mathematics and Statistics, Calgary, Canada}
\author[3,4]{Konstantin Riedl\samethanks[1]}
\affil[3]{Technical University of Munich, School of Computation, Information and Technology, Department of Mathematics, Munich, Germany}
\affil[4]{Munich Center for Machine Learning, Munich, Germany}

\begin{document}
\maketitle
\begin{abstract}
\noindent
	In this paper we propose consensus-based optimization for saddle point problems~(CBO-SP), a novel multi-particle metaheuristic derivative-free optimization method capable of provably finding global Nash equilibria.
	Following the idea of swarm intelligence, the method employs \revised{two groups} of interacting particles, \revised{one of which performs} a minimization over one variable \revised{while the other performs} a maximization over the other variable. \revised{The two groups constantly exchange information through a suitably weighted average.}
	This paradigm permits for a passage to the mean-field limit, which makes the method amenable to theoretical analysis and it allows to obtain rigorous convergence guarantees under reasonable assumptions about the initialization and the objective function, which most notably include nonconvex-nonconcave objectives.
	\revised{We further provide numerical evidence for the success of the algorithm.}
\end{abstract}

{\noindent\small{\textbf{Keywords:} saddle point problems, Nash equilibria, nonconvex-nonconcave, derivative-free optimization, metaheuristics, consensus-based optimization, Fokker-Planck equations}}\\

{\noindent\small{\textbf{AMS subject classifications:} 90C47, 65C35, 65K05, 90C56, 35Q90, 35Q83}}


\section{Introduction} \label{sec:introduction}
Optimization problems where the goal is to find the best possible objective value for the worst-case scenario can be formulated as minimax optimization problems of the form
\begin{equation*}
	\min_{x\in\CX} \max_{y\in\CY} \CE(x,y).
\end{equation*}
To be more specific, given a class of objective functions~$\{\CE(\,\cdot\,,y), y\in\CY\}$, the aim is to determine the argument~$x^*\in\CX$ that leads to the smallest objective value even for the worst-case function parametrized by~$y^*\in\CY$.
Such type of problems were originally formulated in two-player zero-sum game theory~\cite{von2007theory} but now arise in many areas in mathematics, biology, the social sciences and especially in economics~\cite{myerson1997game}. Diverse applications may be found  in engineering, operational research, biology, ecology, finance, economics, energy industry, environmental sciences and so on. In the last few years, minimax optimization has also experienced substantial attention from the signal processing community, due to its connection to distributed processing~\cite{chang2020distributed}, robust transceiver design~\cite{liu2013max}, and communication in the presence of jammers~\cite{gohary2009generalized}.
Moreover, in modern machine learning, several problems are formulated as minimax optimization, such as the training of generative adversarial networks (GANs)~\cite{goodfellow2020generative}, multi-agent reinforcement learning~\cite{omidshafiei2017deep}, fair machine learning~\cite{madras2018learning}, and adversarial training~\cite{madry2018towards}.
For example, when training GANs, $x$ models the parameter of a generator, usually a neural network, whose aim is to generate synthetic data with the same statistics as of a given training set, while $y$ represents the parameters of a competing discriminator, who has to distinguish generated data by the generator from data of the true distribution. 
Relatedly, in adversarial machine learning, one aims at learning the parameters $x$ of a model in a robust manner by exposing it during training to possible adversarial attacks modeled by $y$.
Both examples can be interpreted as a game between two neural networks trained in an adversarial manner until some kind of equilibrium is reached.

In a two-player zero-sum game, the joint payoff function $\CE(x,y)$ encodes the gain of the maximization player whose action is to choose $y\in\CY$, as well as the loss of the minimization player controlling the action $x\in\CX$.
In simultaneous games, each player chooses its action without the knowledge of the action chosen by the other player, so both players act simultaneously.
Conversely, in sequential games there is an intrinsic order according to which the players take their actions, meaning that the ordering of the minimization and maximization \revised{matters, i.e., it plays a priorly a role} whether is $\min_x\max_y$ or $\max_y\min_x$.
GANs and adversarial training, for instance, are in fact sequential games in their standard formulations.
In the classical case, where the payoff function~$\CE$ is \textit{convex-concave} (i.e., $\CE(\,\cdot\,,y)$ is convex for all $y\in\CY$ and $\CE(x,\,\cdot\,)$ is concave for all $x\in\CX$), the intrinsic order of sequential games does not matter under an additional compactness assumption on either $\CX$ or $\CY$ by the well-known minimax theorems of von Neumann and Sion~\cite{v1928theorie,sion1958general}.
However, nowadays, most modern applications in signal processing and machine learning entail the setting of \textit{nonconvex-nonconcave} minimax problems, where the minimization and maximization problems are potentially nonconvex and nonconcave.
This is significantly more complicated and available tool sets and theories are very limited; see the review paper \cite{razaviyayn2020nonconvex}.

A well-known notion of optimality originating from game theory is the one of Nash equilibria (also referred to as saddle points)~\cite{nash1950equilibrium}, where neither of the players has anything to gain by changing only his own strategy.
This concept is formalized within the following definition.
\begin{definition} \label{nash}
	A point $(x^*,y^*)\in\CX\times\CY$ is called Nash equilibrium or saddle point of a function $\CE$ if it holds
	\begin{equation*}
		\CE(x^*,y)\leq \CE(x^*,y^*)\leq \CE(x,y^*) \quad\text{for all } (x,y)\in \CX\times\CY
	\end{equation*}
	or, equivalently, if
	\begin{equation*}
		\min_{x\in\CX} \max_{y\in\CY} \CE(x,y)
		= \CE(\saddlepoint) = 
		\max_{y\in\CY} \min_{x\in\CX} \CE(x,y).
	\end{equation*}
	To keep the notation concise we write $\CE^*$ for $\CE(x^*,y^*)$ in what follows.
\end{definition}

\noindent 
In the \textit{convex-concave} setting an approximate Nash equilibrium can be found efficiently by variants of gradient descent-ascent~(GDA) algorithms~\cite{bubeck2015convex,hazan2016introduction}, which alternate between one or more gradient decent steps in the $x$ variable and gradient ascent steps in the $y$ coordinate. Indeed, even if $\CE(x,y)$ is either concave in $y$ or convex in $x$, there are available some multi-step GDA algorithms; see \cite{nouiehed2019solving,razaviyayn2020nonconvex} for instance. However, as soon as the payoff function becomes \textit{nonconvex-nonconcave}, finding a global equilibrium is in general an NP-hard problem~\cite{murty1987some}.
For this reason, some recent works such as \cite{mazumdar2019finding,daskalakis2018limit} consider a local version of equilibria.
More precisely, a point~$(x^*,y^*)\in\CX\times\CY$ is called local Nash equilibrium if there exists some $\delta>0$ such that $(x^*,y^*)$ satisfies Definition~\ref{nash} in a $\delta$-neighborhood of $(x^*,y^*)$. 
Local Nash equilibria can be characterized in terms of the so-called quasi-Nash equilibrium condition~\cite{pang2011nonconvex} or the first-order Nash equilibrium condition~\cite{nouiehed2019solving}. Even so, we mention two recent works where special classes of nonconvex-nonconcave payoff functions are concerned. When  $\CE(x,y)$  is weakly convex in $x$ and weakly concave in $y$ and the associated Minty variational inequality admits a solution, Liu et al.\@~\cite{liu2021first} employ the inexact proximal point method and prove the first-order convergence, while under the so-called ``sufficiently bilinear" condition,  the stochastic Hamiltonian method is investigated by Loizou et al.\@~\cite{loizou2020stochastic}. In this work, we shall drop such restrictions and the gradient-dependence in the algorithms and consider a \textit{zero-order} (derivative-free) method with rigorous convergence guarantees.
\revised{Note that the family of population-based algorithms, such as Particle Swarm Optimization (PSO)~\cite{kennedy1995particle}, has been adapted to solve min-max problems as done for instance in \cite{shi2002co,krohling2004co,laskari2002particle}.} 
One straightforward approach is to treat the min-max problem as a minimization problem and embed the maximization part in the calculation of the objective values \cite{laskari2002particle}.
Alternatively, a multi-PSO strategy \cite{shi2002co,krohling2004co} may be employed, where the min-max problem is converted into two optimization problems,
one being a maximization problem, and the other a minimization problem.
Two PSO algorithms are then used to solve these two optimization problems, respectively, and they are run independently.
Each PSO is treated as a changing environment of the other PSO, allowing them to cooperate through the calculation of the objective.
Both approaches cannot avoid the necessity of nested loops/circles of optimization algorithms, which significantly increases the time complexity.

In the present paper we propose a zero-order consensus-based optimization method for finding the global Nash equilibrium~$(x^*,y^*)$ of a smooth objective function~$\CE:\CX\times \CY\rightarrow \bbR$ with $\CX=\bbR^{d_1}$ and $\CY=\bbR^{d_2}$, \revised{which is designed to be amenable to a rigorous theoretical convergence analysis, missing so far in the literature on population-based methods for min-max problems.}
The dynamics of the algorithm is inspired by consensus-based optimization, a paradigm for global nonconvex minimizations, which was introduced by the authors of~\cite{pinnau2017consensus}.
Their method employs a system of interacting particles which explore the energy landscape in order to form a global consensus about the global minimizer of the objective function as time passes.
Taking inspiration from this concept, \revised{let us consider two sets of particles~$(X^i)_{i=1}^{N_1}$ and $(Y^i)_{i=1}^{N_2}$ of potentially different size, one for minimization, the other for maximization.
Each individual particle of either set is formally described by a stochastic process.}
In order to achieve consensus about the equilibrium point of $\CE$, the particles interact through a system of stochastic differential equations~(SDEs) of the form
\begin{subequations} \label{eq:saddle_point_dynamics_micro}
\begin{align}
	dX_t^i
		&= -\lambda_1\left(X_t^i - \conspointx{\empmeasureX{t}}\right)dt
		+ \sigma_1 D\!\left(X_t^i - \conspointx{\empmeasureX{t}}\right) dB_t^{X,i},
		&& \textstyle\empmeasureX{t} = \frac{1}{\revised{N_1}}\sum_{i=1}^{\revised{N_1}} \delta_{X_t^i},
		\label{eq:saddle_point_dynamics_micro_X} \\
	dY_t^i
		&= -\lambda_2\left(Y_t^i - \conspointy{\empmeasureY{t}}\right)dt
		+ \sigma_2 D\!\left(Y_t^i - \conspointy{\empmeasureY{t}}\right)  dB_t^{Y,i},
		&& \textstyle\empmeasureY{t} = \frac{1}{\revised{N_2}}\sum_{i=1}^{\revised{N_2}} \delta_{Y_t^i},
	\label{eq:saddle_point_dynamics_micro_Y}
\end{align}
\end{subequations}
which is complemented by suitable initial conditions~\revised{$X_0^i \sim \rho_{X,0}\in\CP(\bbR^{d_1})$ for $i=1,\dots,N_1$ and $Y_0^i \sim \rho_{Y,0}\in\CP(\bbR^{d_2})$ for $i=1,\dots,N_2$} and where $\big(\big(B_t^{X,i}\big)_{t\geq0}\big)_{i=1,\dots,\revised{N_1}}$ and $\big(\big(B_t^{Y,i}\big)_{t\geq0}\big)_{i=1,\dots,\revised{N_2}}$ are independent standard Brownian motions in $\bbR^{d_1}$ and $\bbR^{d_2}$, respectively.
Moreover, $\empmeasureX{t}$ and $\empmeasureY{t}$ denote the empirical measures of the particles' $x$- and $y$-positions, \revised{respectively.}
While the dynamics~\eqref{eq:saddle_point_dynamics_micro_X} performs minimization in the $x$-variable, \eqref{eq:saddle_point_dynamics_micro_Y} performs maximization in the $y$-coordinate.
This is encoded in the computation of the so-called consensus point~$\big(\conspointx{\empmeasureX{t}},\conspointy{\empmeasureY{t}}\big)$, whose components are given by
\begin{subequations} \label{eq:conspoint}
\begin{align} \label{eq:conspoint_x}
	\conspointx{\empmeasureX{t}}
	&= \!\int \! x \,\frac{\omegaa\big(x,\int \!y\,d\empmeasureY{t}(y)\big)}{\Nbig{\omegaa\big(\,\cdot\,,\int \!y\,d\empmeasureY{t}(y)\big)}_{L_1(\empmeasureX{t})}}\,d\empmeasureX{t}(x),
	&& \textrm{with}\ \,
	\omegaa(x,y) \!:=\! \exp(-\alpha \CE(x,y)),\\
	\conspointy{\empmeasureY{t}}
	&= \!\int \! y \,\frac{\omegab\big(\int \!x\,d\empmeasureX{t}(x),y\big)}{\Nbig{\omegab\big(\int \!x\,d\empmeasureX{t}(x),\,\cdot\,\big)}_{L_1(\empmeasureY{t})}}\,d\empmeasureY{t}(y),
	&& \textrm{with}\ \,
	\omegab(x,y) \!:=\! \exp(\beta \CE(x,y)).\label{eq:conspoint_y}
\end{align}
\end{subequations}
Attributed to the Laplace principle~\cite{miller2006applied}, $\conspointx{\empmeasureX{t}}$ can be interpreted as an approximation of $\argmin_{i=1,\dots,\revised{N_1}} \CE(X_t^i, \int \!y\,d\empmeasureY{t}(y))$ as $\alpha\rightarrow\infty$ while $\conspointy{\empmeasureY{t}}\approx\argmax_{i=1,\dots,\revised{N_2}} \CE(\int \!x\,d\empmeasureX{t}(x),Y_t^i)$ as $\beta\rightarrow\infty$, see, e.g., \cite[Equation~(7)]{fornasier2021consensus}.
The dynamics of each of the particles in~\eqref{eq:saddle_point_dynamics_micro} is governed by two terms.
A drift term drags the particles towards \revised{the respective component of} the instantaneous consensus point~$\big(\conspointx{\empmeasureX{t}},\conspointy{\empmeasureY{t}}\big)$ and thereby expectedly improves the position of the particles.
The second term injects stochasticity into the dynamics by diffusing the particles according to a scaled Brownian motion, which features the exploration of the landscape of the objective.
In what follows we use anisotropic noise, i.e., $D(\,\cdot\,) = \diag(\,\cdot\,)$, \revised{which is typically more competitive in high dimensions compared to isotropic noise~$D(\,\cdot\,) = \N{\,\cdot\,}_2$, see, e.g., \cite{carrillo2019consensus,fornasier2021convergence}.
The theoretical results of this paper, however, can be obtained mutatis mutandis also in the isotropic setting.}

\revised{An implementable scheme for a numerical algorithm can be obtained from \eqref{eq:saddle_point_dynamics_micro} by a simple Euler-Maruyama time discretization~\cite{higham2001algorithmic,platen1999introduction}.
For details about the implementation we refer to Algorithm~\ref{algorithm:CBOSP} in Section~\ref{sec:implementation}.}

\revised{\begin{remark}
	While the definition of the consensus point in~\eqref{eq:conspoint} is a natural option, there are two equally reasonable alternatives.
	The first possibility is to replace the mean~$\int \!y\,d\empmeasureY{t}(y)$ in~\eqref{eq:conspoint_x} simply by $y$ and integrate w.r.t.\@ the joint measure $\empmeasure{t}$. This case would require $N_1=N_2$.
	Analogously, $\int \!x\,d\empmeasureX{t}(x)$ is substituted by $x$ in \eqref{eq:conspoint_y}.
	The second option is to use the other component of the consensus point instead of the respective mean, i.e., $\conspointy{\empmeasureY{t}}$ replaces $\int \!y\,d\empmeasureY{t}(y)$ in~\eqref{eq:conspoint_x} and $\conspointx{\empmeasureX{t}}$ substitutes $\int \!x\,d\empmeasureX{t}(x)$ in~\eqref{eq:conspoint_y}.
	
	The main motivation for using the variant as in~\eqref{eq:conspoint} is of theoretical nature.
	Using either of the other two alternatives significantly complicates the convergence analysis in Sections~\ref{sec:convergence} and \ref{sec:proofs}.
\end{remark}}

Understanding the convergence properties of the dynamics~\eqref{eq:saddle_point_dynamics_micro} can take place either by investigating the long time behavior of the interacting particle system itself, or by analyzing the macroscopic behavior of the agent density associated with~\eqref{eq:saddle_point_dynamics_micro} through a mean-field limit.
This theoretical approach proved successful in \cite{carrillo2018analytical,fornasier2021consensus,carrillo2019consensus,fornasier2021convergence,riedl2022leveraging,fornasier2020consensus_sphere_convergence,fornasier2021anisotropic,carrillo2021consensus,borghi2021constrained,borghi2022consensus} for proving global convergence for several variants of consensus-based optimization in the setting of minimization.
It is moreover theoretically justified by the mean-field approximation which shows that $(\empmeasureX{t},\empmeasureY{t})$ converges in some sense to a mean field law~$(\rho_{X,t},\rho_{Y,t})$ as $\revised{N_1,N_2}\rightarrow\infty$.
Again, for consensus-based optimization there exist by now several results in this direction such as \cite{fornasier2021consensus,huang2021MFLCBO,fornasier2020consensus_hypersurface_wellposedness}, which may be extended to CBO-SP in an immediate manner.
In the setting of saddle point problems, the mean-field dynamics associated with~\eqref{eq:saddle_point_dynamics_micro}  can be described by the self-consistent mono-particle dynamics
\begin{subequations} \label{eq:saddle_point_dynamics_macro}
\begin{align}
	d\OX_t
		&= -\lambda_1\left(\OX_t - \conspointx{\rho_{X,t}}\right)dt
		+ \sigma_1 D\!\left(\OX_t - \conspointx{\rho_{X,t}}\right) dB_t^{X}, 
		&& \textstyle\rho_{X,t} = \int d\rho_{t}(\,\cdot\,,y),
		\label{eq:saddle_point_dynamics_macro_X} \\
	d\OY_t
		&= -\lambda_2\left(\OY_t - \conspointy{\rho_{Y,t}}\right)dt
		+ \sigma_2 D\!\left(\OY_t - \conspointy{\rho_{Y,t}}\right)  dB_t^{Y},
		&& \textstyle\rho_{Y,t} = \int d\rho_{t}(x,\,\cdot\,),
	\label{eq:saddle_point_dynamics_macro_Y}
\end{align}
\end{subequations}
where $\rho_t = \rho(t)= \Law\left((\OX_t,\OY_t)\right)$ with marginals $\rho_{X,t}$ and $\rho_{Y,t}$, respectively.
In particular, the measure $\rho\in\CC([0,T],\CP(\bbR^{d_1+d_2}))$ weakly satisfies the nonlinear nonlocal Fokker-Planck equation
\begin{equation} \label{eq:fokker_planck}
\begin{split}
	\partial_t\rho_t
	= &\, \lambda_1\divergence_x \big(\!\left(x - \conspointx{\rho_t^X}\right)\rho_t\big)
		+ \lambda_2\divergence_y \big(\!\left(y - \revised{\conspointy{\rho_t^Y}}\right)\rho_t\big) \\
	& + \frac{\sigma_1^2}{2}\sum_{k=1}^{d_1} \partial_{x_kx_k}^2\big((x - \conspointx{\rho_t^X})_k^2 \rho_t\big)
		+ \frac{\sigma_2^2}{2}\sum_{k=1}^{d_2} \partial_{y_ky_k}^2\big((y - \conspointy{\rho_t^Y})_k^2 \rho_t\big).
\end{split}
\end{equation}

\paragraph{Contributions.} 

Motivated by the fundamental importance of \textit{nonconvex-nonconcave} saddle point problems in various applicational areas and the desire for having numerical algorithms with rigorous global convergence guarantees, we theoretically analyze in this work a novel consensus-based optimization method (CBO-SP) capable of tackling saddle point problems.
Using mean-field analysis techniques, we rigorously prove that CBO-SP converges to saddle points as the number of interacting particles goes to infinity.
Our results hold under reasonable assumptions about the objective function and under certain conditions of well-preparation of the hyperparameters and the initial data.

\subsection{Organization}

In Section~\ref{sec:well-posendness} we first investigate the well-posedness of both the interacting particle system~\eqref{eq:saddle_point_dynamics_micro} of CBO-SP and its associated mean-field dynamics~\eqref{eq:saddle_point_dynamics_macro}.
Section~\ref{sec:convergence} then presents and discusses the main theoretical statement of this work concerned with the convergence of the mean-field dynamics~\eqref{eq:saddle_point_dynamics_macro} towards saddle points of the objective function~$\CE$, which are proven in Section~\ref{sec:proofs}.
Section~\ref{sec:numerics} \revised{contains details about the implementation of the numerical algorithm as well as} instructive numerical examples which illustrate how CBO-SP works, before we conclude the paper in Section~\ref{sec:conclusion}.
In the GitHub repository \url{https://github.com/KonstantinRiedl/CBOSaddlePoints} we provide the Matlab code implementing CBO-SP.

\section{Well-Posedness of CBO-SP and its Mean-Field Dynamics} \label{sec:well-posendness}

In the first part of this section we provide a well-posedness result about the interacting particle system~\eqref{eq:saddle_point_dynamics_micro} of CBO-SP, i.e., we show that a process obeying \eqref{eq:saddle_point_dynamics_micro} exists and is unique.
Afterwards we also prove the well-posedness of the nonlinear macroscopic SDE~\eqref{eq:saddle_point_dynamics_macro}.

\subsection{Well-Posedness of the Interacting Particle System} \label{subsec:well-posendness_interacting_particle_system}

To keep the notation concise in what follows, let us denote the state vector of the entire particle system~\eqref{eq:saddle_point_dynamics_micro} by $\mathbf{Z}\in\CC([0,\infty),\bbR^{\revised{N_1d_1+N_2d_2}})$ with $\mathbf{Z}(t) = \mathbf{Z}_t = \left((X_t^1)^T, \dots, (X_t^{\revised{N_1}})^T, (Y_t^1)^T, \dots, (Y_t^{\revised{N_2}})^T\right)^T$ for every $t\geq0$.
Equation~\eqref{eq:saddle_point_dynamics_micro} can then be reformulated as
\begin{equation} \label{eq:saddle_point_dynamics_micro_compact}
	d\mathbf{Z}_t
	= -\bm{\lambda}\mathbf{F}(\mathbf{Z}_t)\,dt + \bm{\sigma}\mathbf{M}(\mathbf{Z}_t)\,d\mathbf{B}_t
\end{equation}
with $(\mathbf{B}_t)_{t\geq0}$ being a standard Brownian motion in $\bbR^{\revised{N_1d_1+N_2d_2}}$ and  definitions
\begin{equation*}
\begin{split}
	\mathbf{F}(\mathbf{Z}_t)
	:= \,&\left({F}^{1,X}(\mathbf{Z}_t)^T, \dots, {F}^{\revised{N_1},X}(\mathbf{Z}_t)^T,{F}^{1,Y}(\mathbf{Z}_t)^T, \dots, {F}^{\revised{N_2},Y}(\mathbf{Z}_t)^T\right)^T \\
	&\text{ with } {F}^{i,X}(\mathbf{Z}_t) = \left(X_t^i - \conspointx{\empmeasureX{t}}\right)\! \text{ and } {F}^{i,Y}(\mathbf{Z}_t) = \left(Y_t^i - \conspointy{\empmeasureY{t}}\right), \\
	\mathbf{M}(\mathbf{Z}_t)
	:= \,&\diag\left({M}^{1,X}(\mathbf{Z}_t), \dots, {M}^{\revised{N_1},X}(\mathbf{Z}_t), {M}^{1,Y}(\mathbf{Z}_t), \dots, {M}^{\revised{N_2},Y}(\mathbf{Z}_t)\right) \\
	&\text{ with } {M}^{i,X}(\mathbf{Z}_t) = D\!\left(X_t^i - \conspointx{\empmeasureX{t}}\right)\! \text{ and } {M}^{i,Y}(\mathbf{Z}_t) = D\!\left(Y_t^i - \conspointy{\empmeasureY{t}}\right).
\end{split}
\end{equation*}
The $\diag$ operator in the definition of $\mathbf{M}$ maps the input matrices onto a block-diagonal matrix with them as its diagonal.
$\bm{\lambda}$ and $\bm{\sigma}$ are $\revised{(N_1d_1+N_2d_2) \times (N_1d_1+N_2d_2)}$-dimensional diagonal matrices, whose first $\revised{N_1}d_1$ entries are $\lambda_1$ and $\sigma_1$, and the remaining $\revised{N_2}d_2$ entires are $\lambda_2$ and $\sigma_2$, respectively.

Having fixed the notation, we have the following well-posedness result for the SDE system~\eqref{eq:saddle_point_dynamics_micro_compact}, respectively \eqref{eq:saddle_point_dynamics_micro}, which is proven towards the end of this section.

\begin{theorem} \label{theorem:well-posedness_interacting particle dynamcis}
	Let $\CE\in\CC(\bbR^{d_1+d_2})$ be locally Lipschitz continuous.
	Then, for $\revised{N_1},\revised{N_2}\in\bbN$ fixed, the system of SDEs~\eqref{eq:saddle_point_dynamics_micro} admits a unique strong solution~$\left(\mathbf{Z}_t\right)_{t\geq0}$ for any initial condition~$\mathbf{Z}_0$ satisfying $\bbE\Nnormal{\mathbf{Z}_0}_2^2<\infty$.
\end{theorem}

In order to employ the standard result~\cite[Chapter~5, Theorem~3.1]{Durrett} about the existence and uniqueness of solutions to SDEs, we need to verify that the coefficients of the SDE are locally Lipschitz continuous and of at most linear growth.
This is inherited from the assumed local Lipschitz continuity of $\CE$ as we make explicit in the subsequent lemma.

\begin{lemma} \label{theorem:linear growth well-posedness_interacting particle dynamcis}
	Let $\revised{N_1,N_2}\in\bbN$, $\alpha,\beta>0$ and $R>0$ be arbitrary.
	Let $\mathbf{z},\mathbf{\widehat{z}}\in\bbR^{\revised{N_1d_1+N_2d_2}}$ be of the form $\mathbf{z} = (\mathbf{x}^T,\mathbf{y}^T)^T = \big((x^1)^T,\dots,(x^{\revised{N_1}})^T,(y^1)^T,\dots,(y^{\revised{N_2}})^T\big)^T$ and analogously for $\mathbf{\widehat{z}}$.
	Then, for any $\mathbf{z},\mathbf{\widehat{z}}$ with $\Nnormal{\mathbf{z}}_2 \leq R$ and $\Nnormal{\mathbf{\widehat{z}}}_2 \leq R$ it \revised{holds} \revised{for any $i$} the bounds
	\begin{equation*}
		\Nbig{{F}^{i,X}(\mathbf{z})}_2
		\leq \Nnormal{x^i}_2 + \Nnormal{\mathbf{x}}_2
		\quad\text{and}\quad
		\Nbig{{F}^{i,Y}(\mathbf{z})}_2
		\leq \Nnormal{y^i}_2 + \Nnormal{\mathbf{y}}_2,
	\end{equation*}
	and, abbreviating $c_R(\gamma) := 4\gamma e^{2\gamma\Delta_R\CE}\Nbig{\N{\nabla_\mathbf{z}\CE}_2}_{L^\infty(B_R)}$
	with $\Delta_R\CE := \sup_{\mathbf{z}\in B_R}\CE(\mathbf{z})-\inf_{\mathbf{z}\in B_R}\CE(\mathbf{z})$,
	\begin{equation*}
	\begin{split}
		\Nbig{{F}^{i,X}(\mathbf{z})-{F}^{i,X}(\mathbf{\widehat{z}})}_2
		&\leq \Nnormal{x^i-\widehat{x}^i}_2 + \left(1+\frac{c_R(\alpha)}{\revised{N_1}} \sqrt{\revised{N_1}\N{\widehat{x}^i}_2^2+\N{\mathbf{\widehat{x}}}_2^2}\right) \big(\!\N{\mathbf{x}-\mathbf{\widehat{x}}}_2 + \N{\mathbf{y}-\mathbf{\widehat{y}}}_2\big), \\
		\Nbig{{F}^{i,Y}(\mathbf{\mathbf{z}})-{F}^{i,Y}(\mathbf{\widehat{z}})}_2
		& \leq \Nnormal{y^i-\widehat{y}^i}_2 + \left(1+\frac{c_R(\beta)}{\revised{N_2}} \sqrt{\revised{N_2}\N{\widehat{y}^i}_2^2+\N{\mathbf{\widehat{y}}}_2^2}\right) \big(\!\N{\mathbf{x}-\mathbf{\widehat{x}}}_2 + \N{\mathbf{y}-\mathbf{\widehat{y}}}_2\big).
	\end{split}
	\end{equation*}
\end{lemma}

\begin{proof}
	To derive the first bound we note that
	\begin{equation*}
		\Nbig{{F}^{i,X}(\mathbf{z})}_2
		= \N{x^i - \sum_{j=1}^{\revised{N_1}} x^j \frac{\omegaa\big(x^j, \frac{1}{\revised{N_1}} \sum_{k=1}^{\revised{N_1}} y^k\big)}{\sum_{j=1}^{\revised{N_1}} \omegaa\big(x^j, \frac{1}{\revised{N_1}} \sum_{k=1}^{\revised{N_1}} y^k\big)}}_2
		\leq \Nnormal{x^i}_2 + \N{\mathbf{x}}_2.
	\end{equation*}
	Analogously the bound for $\Nbig{{F}^{i,Y}(\mathbf{z})}_2$ is obtained.
	For the other estimates we first notice that
	\begin{align*}
		{F}^{i,X}(\mathbf{z})-{F}^{i,X}(\mathbf{\widehat{z}})
		&= \frac{\sum_{j=1}^{\revised{N_1}} (x^i-x^j) \, \omegaa\big(x^j, \frac{1}{\revised{N_1}} \sum_{k=1}^{\revised{N_1}} y^k\big)}{\sum_{j=1}^{\revised{N_1}}\omegaa\big(x^j, \frac{1}{\revised{N_1}} \sum_{k=1}^{\revised{N_1}} y^k\big)}
			- \frac{\sum_{j=1}^{\revised{N_1}} (\widehat{x}^i-\widehat{x}^j) \,\omegaa\big(\widehat{x}^j, \frac{1}{\revised{N_1}} \sum_{k=1}^{\revised{N_1}} \widehat{y}^k\big)}{\sum_{j=1}^{\revised{N_1}}\omegaa\big(\widehat{x}^j, \frac{1}{\revised{N_1}} \sum_{k=1}^{\revised{N_1}} \widehat{y}^k\big)} \\
		&= I_1 + I_2 + I_3
	\end{align*}
	with $I_1$, $I_2$ and $I_3$ being defined as in what follows.
	Firstly, for $I_1$ we have
	\begin{align*}
		\N{I_1}_2 &:= \N{\frac{\sum_{j=1}^{\revised{N_1}} \left((x^i-x^j) - (\widehat{x}^i-\widehat{x}^j)\right) \omegaa\big(x^j, \frac{1}{\revised{N_1}} \sum_{k=1}^{\revised{N_1}} y^k\big)}{\sum_{j=1}^{\revised{N_1}} \omegaa\big(x^j, \frac{1}{\revised{N_1}} \sum_{k=1}^{\revised{N_1}} y^k\big)}}_2
		\leq \Nnormal{x^i-\widehat{x}^i}_2 + \N{\mathbf{x}-\mathbf{\widehat{x}}}_2.
	\end{align*}
	For $I_2$ and $I_3$, on the other hand, let us first notice that it \revised{holds}
	\begin{align*}
		&\abs{\omegaa\left(x^j, \frac{1}{\revised{N_1}} \sum_{k=1}^{\revised{N_1}} y^k\right) - \omegaa\left(\widehat{x}^j, \frac{1}{\revised{N_1}} \sum_{k=1}^{\revised{N_1}} \widehat{y}^k\right)} \\
		&\qquad\quad\leq \alpha e^{-\alpha\inf_{\mathbf{z}\in B_R}\CE(\mathbf{x},\mathbf{y})}\Nbig{\N{\nabla_\mathbf{z}\CE}_2}_{L^\infty(B_R)} \left(\N{x^j-\widehat{x}^j}_2 + \frac{1}{\revised{N_1}} \sum_{k=1}^{\revised{N_1}}\Nbig{y^k-\widehat{y}^k}_2 \right)
	\end{align*}
	and
	\begin{equation*}
		\frac{1}{\sum_{j=1}^{\revised{N_1}}\omegaa\big(x^j, \frac{1}{\revised{N_1}} \sum_{k=1}^{\revised{N_1}} y^k\big)}
		\leq \frac{1}{\revised{N_1}\inf_{\mathbf{z}\in B_R}\exp(-\alpha\CE(\mathbf{x},\mathbf{y}))}
		\leq \frac{1}{\revised{N_1} e^{-\alpha\sup_{\mathbf{z}\in B_R}\CE(\mathbf{x},\mathbf{y})}}.
	\end{equation*}
	With this we immediately obtain for the norm of $I_2$ the upper bound
	\begin{align*}
		\N{I_2}_2 &:= \N{\frac{\sum_{j=1}^{\revised{N_1}} (\widehat{x}^i-\widehat{x}^j) \left(\omegaa\big(x^j, \frac{1}{\revised{N_1}} \sum_{k=1}^{\revised{N_1}} y^k\big) - \omegaa\big(\widehat{x}^j, \frac{1}{\revised{N_1}} \sum_{k=1}^{\revised{N_1}} \widehat{y}^k\big)\right)}{\sum_{j=1}^{\revised{N_1}}\omegaa\big(x^j, \frac{1}{\revised{N_1}} \sum_{k=1}^{\revised{N_1}} y^k\big)}}_2\\
		&\leq \frac{2\alpha e^{\alpha\Delta_R\CE}\Nbig{\N{\nabla_\mathbf{z}\CE}_2}_{L^\infty(B_R)}}{\revised{N_1}} \sqrt{\revised{N_1}\N{\widehat{x}^i}_2^2+\N{\mathbf{\widehat{x}}}_2^2} \big(\!\N{\mathbf{x}-\mathbf{\widehat{x}}}_2 + \N{\mathbf{y}-\mathbf{\widehat{y}}}_2\!\big),
	\end{align*}
	where we abbreviate $\Delta_R\CE := \sup_{\mathbf{z}\in B_R}\CE(\mathbf{x},\mathbf{y})-\inf_{\mathbf{z}\in B_R}\CE(\mathbf{x},\mathbf{y})$.
	Similarly, for $I_3$ we have
	\begin{align*}
		\N{I_3}_2 &:= \N{\sum_{j=1}^{\revised{N_1}} (\widehat{x}^i-\widehat{x}^j) \, \omegaa\!\left(\widehat{x}^j, \frac{1}{\revised{N_1}} \sum_{k=1}^{\revised{N_1}} \widehat{y}^k\right) \frac{\left(\sum_{j=1}^{\revised{N_1}} \omegaa\big(\widehat{x}^j, \frac{1}{\revised{N_1}} \sum_{k=1}^{\revised{N_1}} \widehat{y}^k\big) - \omegaa\big(x^j, \frac{1}{\revised{N_1}} \sum_{k=1}^{\revised{N_1}} y^k\big)\right)}{\sum_{j=1}^{\revised{N_1}} \omegaa\big(\widehat{x}^j, \frac{1}{\revised{N_1}} \sum_{k=1}^{\revised{N_1}} \widehat{y}^k\big)\sum_{j=1}^{\revised{N_1}} \omegaa\big(x^j, \frac{1}{\revised{N_1}} \sum_{k=1}^{\revised{N_1}} y^k\big)}}_2\\
		&\leq \frac{2\alpha e^{2\alpha\Delta_R\CE}\Nbig{\N{\nabla_\mathbf{z}\CE}_2}_{L^\infty(B_R)}}{\revised{N_1}} \sqrt{\revised{N_1}\N{\widehat{x}^i}_2^2+\N{\mathbf{\widehat{x}}}_2^2} \big(\!\N{\mathbf{x}-\mathbf{\widehat{x}}}_2 + \N{\mathbf{y}-\mathbf{\widehat{y}}}_2\!\big).
	\end{align*}
	Combining these bounds yields the result.
	Analogously we can bound $\Nbig{{F}^{i,Y}(\mathbf{\mathbf{z}})-{F}^{i,Y}(\mathbf{\widehat{z}})}_2$.
\end{proof}


\begin{proof}[Proof of Theorem~\ref{theorem:well-posedness_interacting particle dynamcis}]
	The statement follows by invoking the standard result~\cite[Chapter~5, Theorem~3.1]{Durrett} \revised{(see Theorem~\ref{thm31} in the Appendix)} on the existence and pathwise uniqueness of a strong solution.
	That \revised{Condition~\!\ref{thm31_i} of Theorem~\ref{thm31}} about the local Lipschitz continuity and linear growth of $\mathbf{F}(\mathbf{Z}_t)$ and $\mathbf{M}(\mathbf{Z}_t)$ holds, follows immediately from Lemma~\ref{theorem:linear growth well-posedness_interacting particle dynamcis}.
	To ensure \revised{Condition~\!\ref{thm31_ii} of Theorem~\ref{thm31}} we make use of \cite[Chapter~5, Theorem~3.2]{Durrett} \revised{(see Theorem~\ref{thm32} in the Appendix)}  and verify that there exists a constant $b_{\revised{N_1},\revised{N_2}}>0$ such that $-2\bm{\lambda}\mathbf{Z}_t\cdot\mathbf{F}(\mathbf{Z}_t) + \mathrm{tr}(\bm{\sigma}\mathbf{M}(\mathbf{Z}_t)\mathbf{M}(\mathbf{Z}_t)^T\bm{\sigma}^T) \leq b_{\revised{N_1},\revised{N_2}} (1 + \N{\mathbf{Z}_t}_2^2)$.
	Indeed, since 
	\revised{\begin{align*}
		-\bm{\lambda}\mathbf{Z}_t \cdot \mathbf{F}(\mathbf{Z}_t)
		&\leq \lambda_1\sum_{i=1}^{N_1} \N{X_t^i}_2\N{{F}^{i,X }(\mathbf{Z}_t)}_2  + \lambda_2\sum_{i=1}^{N_2} \N{Y_t^i}_2\N{{F}^{i,Y }(\mathbf{Z}_t)}_2  \\
		&\leq \left(\lambda_1\big(1 + \sqrt{N_1}\big) + \lambda_2\big(1 + \sqrt{N_2}\big)\right)  \N{\mathbf{Z}_t}_2^2
	\end{align*}}
	and
	\revised{\begin{align*}
		\mathrm{tr}(\bm{\sigma}\mathbf{M}(\mathbf{Z}_t)\mathbf{M}(\mathbf{Z}_t)^T\bm{\sigma}^T)
		&= \sigma_1^2 \sum_{i=1}^{N_1} \N{{F}^{i,X }(\mathbf{Z}_t)}_2^2  +  \sigma_2^2 \sum_{i=1}^{N_2}  \N{{F}^{i,Y }(\mathbf{Z}_t)}_2^2  \\
		&\leq 2\left(\sigma_1^2\big(1 + {N_1}\big) + \sigma_2^2\big(1 + {N_2}\big)\right)  \N{\mathbf{Z}_t}_2^2,
	\end{align*}}
	the former holds with $b_{\revised{N_1},\revised{N_2}}$ defined as the sum of the two former upper bounds.
\end{proof}

\subsection{Well-Posedness of the Mean-Field Dynamics} \label{subsec:well-posendness_mean_field_limit}

In what follows let us furthermore ensure the well-posedness of the mean-field dynamics~\eqref{eq:saddle_point_dynamics_macro} and~\eqref{eq:fokker_planck}, which is the main object of our studies in Section~\ref{sec:convergence}.
We prove existence and uniqueness of a solution for objective functions~$\CE$ that satisfies the following conditions.

\begin{definition}[Assumptions] \label{ass:sec:convergence}
	In this section we consider functions $\CE\in\CC^1(\bbR^{d_1+d_2})$, which
	\begin{enumerate}[label=W\arabic*,labelsep=10pt,leftmargin=35pt]
		\item \label{ass:sec:wp_1}
			are bounded in the sense that there exist~$\underbar{\CE}\in\CC^1(\bbR^{d_2})$ and $\overbar{\CE}\in\CC^1(\bbR^{d_1})$ such that
			\begin{equation*}
				\underbar{\CE}(y)\leq\CE(x,y)\leq\overbar{\CE}(x) \quad\text{for all } (x,y)\in\bbR^{d_1+d_2}.
			\end{equation*}
		\item \label{ass:sec:wp_2}
			are locally Lipschitz continuous \revised{in the sense that} there exists a constant $C_1>0$ such that for all $(x,y),(x',y')\in\bbR^{d_1+d_2}$ it holds
			\begin{equation*}
				\SN{\CE(x,y)-\CE(x',y')}
				\leq C_1\!\left(1+\N{x}_2 \!+\! \Nnormal{x'}_2 \!+\! \N{y}_2 \!+\! \Nnormal{y'}_2\right)\left(\Nnormal{x-x'}_2 + \Nnormal{y-y'}_2\right).
			\end{equation*}
		\item \label{ass:sec:wp_3}
			have at most quadratic growth in the sense that there exists a constant $C_2>0$ obeying 
			\begin{subequations}
			\begin{equation*}
			\begin{split}
				\CE(x,y) \!-\! \underbar\CE(y\!+\!sy') &\leq C_2 \big(1\!+\!\N{x}_2^2 \!+\! \N{y}_2^2 \!+\! \Nnormal{y'}_2^2 \big) \quad\text{for all } (x,y),(x',y')\in\bbR^{d_1+d_2}, s\!\in\![0,1],\\
				\overbar\CE(x\!+\!sx') \!-\! \CE(x,y) &\leq C_2 \big(1\!+\!\N{x}_2^2 \!+\! \Nnormal{x'}_2^2 \!+\! \N{y}_2^2\!\big)\quad\text{for all } (x,y),(x',y')\in\bbR^{d_1+d_2}, s\!\in\![0,1].
			\end{split}
			\end{equation*}
			\end{subequations}
	\end{enumerate}
\end{definition}

For such objectives we have the following well-posedness result for the macroscopic SDE~\eqref{eq:saddle_point_dynamics_macro} and its associated Fokker-Planck equation~\eqref{eq:fokker_planck}.

\begin{theorem} \label{theorem:well-posedness_mean-field}
	Let $\CE$ satisfy Assumptions~\ref{ass:sec:wp_1}--\ref{ass:sec:wp_3}.
	Let $T > 0$, $\rho_0 \in \CP_4(\bbR^{d_1+d_2})$.
	Then there exists a unique nonlinear process $(\overbar{X},\overbar{Y}) \in \CC([0,T],\bbR^{d_1+d_2})$ satisfying \eqref{eq:saddle_point_dynamics_macro}.
	The associated law $\rho = \Law(\overbar{X},\overbar{Y})$ has regularity $\rho \in \CC([0,T], \CP_4(\bbR^{d_1+d_2}))$ and is a weak solution to the Fokker-Planck equation~\eqref{eq:fokker_planck}.
\end{theorem}

Before giving the proof of Theorem~\ref{theorem:well-posedness_mean-field} let us first provide some auxiliary results.

\begin{lemma} \label{lem:aux_lemma_bound}
	Let $\revised{\varrho,\widehat\varrho}\in\CP_2(\bbR^{d_1+d_2})$ with $\iint \N{x}_2^2 + \N{y}_2^2 d\revised{\varrho}(x,y) \leq K$ and $\iint \Nnormal{\widehat{x}}_2^2 + \Nnormal{\widehat{y}}_2^2 \,d\revised{\widehat\varrho}(\widehat{x},\widehat{y})  \leq K$.
	Then, under \revised{Assumptions~\ref{ass:sec:wp_1} and~\ref{ass:sec:wp_3}} on~$\CE$, it \revised{holds} for any $s\in[0,1]$ that
	\begin{subequations}
	\begin{align}
		\frac{\exp\left(-\alpha\underbar\CE\left(\int y d\varrho_Y(y)+ s\big(\!\int \widehat y d\revised{\widehat\varrho_Y}(\revised{\widehat y})-\int y d\varrho_Y(y)\big)\right)\right)}{\int \!\omegaa\big(x,\int y d\varrho_Y(y)\big)\,d\varrho_X(x)}
		&\leq \exp\left(\alpha C_2(1+2K)\right):=C_K^\alpha \label{eq:aux_lemma_bound_1} \\
		\intertext{and}
		\frac{\exp\left(\beta\overbar\CE\left(\int x d\varrho_X(x) + s\big(\!\int \widehat x d\revised{\widehat\varrho_X}(\widehat x)-\int x d\varrho_X(x)\big)\right)\right)}{\int \!\omegab\big(\!\int x d\varrho_X(x),y\big)\,d\varrho_X(x)}
		&\leq\exp\left(\beta C_2(1+2K)\right):= C_K^\beta \label{eq:aux_lemma_bound_2}.
	\end{align}
	\end{subequations}
\end{lemma}

\begin{proof}
	By exploiting Assumption~\ref{ass:sec:wp_3} and utilizing Jensen's inequality, we obtain
	\begin{align*}
		&\frac{\exp\left(-\alpha\underbar\CE\left(\int y d\varrho_Y(y) + s\big(\!\int \widehat y d\revised{\widehat\varrho_Y}(\widehat y)-\int y d\varrho_Y(y)\big)\right)\right)}{\int \!\omegaa\big(x,\int y d\varrho_Y(y)\big)\,d\varrho_X(x)} \\
		&\qquad\leq \frac{1}{\exp\left(-\alpha C_2 \left( \int 1+\N{x}_2^2 + \Nbig{\int y d\varrho_Y(y)}_2^2 + \Nbig{\int \widehat y d\revised{\widehat\varrho_Y}(\widehat y)}_2^2 \,d\varrho_X(x)\right)\right)} 
		\leq \exp\left(\alpha C_2(1+2K)\right),
	\end{align*}
	where in the last inequality we integrated the moment bounds on $\varrho$ and $\widehat\varrho$.
	A similar estimate gives~\eqref{eq:aux_lemma_bound_2} by exploiting Assumption~\ref{ass:sec:wp_3}.
\end{proof}

\begin{lemma}
	\label{lem:stabiliy_conspoint}
	Let $\CE$ satisfy Assumption~\ref{ass:sec:wp_1} and the assumptions of Theorem~\ref{theorem:well-posedness_mean-field}.
	Let $\revised{\varrho,\widehat\varrho}\in\CP_4(\bbR^{d_1+d_2})$ with $\iint \N{x}_2^4 \!+\! \N{y}_2^4 d\revised{\varrho}(x,y) \!\leq\! K$ and $\iint \Nnormal{\widehat{x}}_2^4 \!+\! \Nnormal{\widehat{y}}_2^4 \,d\revised{\widehat\varrho}(\widehat{x},\widehat{y})  \!\leq\! K$.
	Then it holds the stability estimate
	\begin{align}
		\N{\conspointx{\varrho_X}-\conspointx{\widehat\varrho_X}}_2 + \N{\conspointy{\varrho_Y}-\conspointy{\widehat\varrho_Y}}_2
		\leq c_0 W_2(\varrho,\widehat\varrho)
	\end{align}
	with $c_0$ depending only on $\alpha,\beta,C_1,C_2$ and $K$.
\end{lemma}

\begin{proof}
		To keep the notation concise we write $\mathsf{E}\varrho_Y:= \int y \,d\varrho_Y(y)$ and $\mathsf{E}\widehat\varrho_Y:= \int \widehat{y} \,d\widehat\varrho_Y(\widehat{y})$ in what follows.
	According to the definition of the consensus point in~\eqref{eq:conspoint_x} we have 
	\begin{align*}
		\conspointx{\varrho_X}-\conspointx{\widehat\varrho_X}
		&= \iint \underbrace{\frac{x \omegaa\big(x,\mathsf{E}\varrho_Y\big)}{\int \!\omegaa\big(x,\mathsf{E}\varrho_Y\big)\,d\varrho_X(x)} - \frac{\widehat{x} \omegaa\big(\widehat{x},\mathsf{E}\widehat\varrho_Y\big)}{\int \!\omegaa\big(\widehat{x},\mathsf{E}\widehat\varrho_Y\big)\,d\widehat\varrho_X(\widehat{x})}}_{=:\, h(x)-h(\widehat{x})} d\pi(x,y,\widehat{x},\widehat{y}),
	\end{align*}
	where $\pi\in\Pi(\varrho,\widehat\varrho)$ is any coupling of $\varrho$ and $\widehat\varrho$.
	By adding and subtracting mixed terms, we obtain the decomposition
	\begin{align*}
		h(x)-h(\widehat{x})
		&= \frac{(x-\widehat{x})\, \omegaa\big(x,\mathsf{E}\varrho_Y\big)}{\int \!\omegaa\big(x,\mathsf{E}\varrho_Y\big)\,d\varrho_X(x)}
			+ \frac{\widehat{x} \big(\omegaa\big(x,\mathsf{E}\varrho_Y\big)-\omegaa\big(\widehat{x},\mathsf{E}\widehat\varrho_Y\big)\big)}{\int \!\omegaa\big(x,\mathsf{E}\varrho_Y\big)\,d\varrho_X(x)} \\
			&\quad\, + \frac{\iint \!\omegaa\big(\widehat{x},\mathsf{E}\widehat\varrho_Y\big)-\omegaa\big(x,\mathsf{E}\varrho_Y\big)\,d\pi(x,y,\widehat{x},\widehat{y})}{\left(\int \!\omegaa\big(x,\mathsf{E}\varrho_Y\big)\,d\varrho_X(x)\right)\left(\int \!\omegaa\big(\widehat{x},\mathsf{E}\widehat\varrho_Y\big)\,d\widehat\varrho_X(\widehat{x})\right)} \widehat{x} \omegaa\big(\widehat{x},\mathsf{E}\widehat\varrho_Y\big)
			=: I_1 + I_2 + I_3,
	\end{align*}
	where $I_1,I_2$ and $I_3$ correspond to the three summands.
	For $I_1$, by using Assumption~\ref{ass:sec:wp_3} and Lemma~\ref{lem:aux_lemma_bound} with $s=0$, we obtain
	\begin{align*}
		\N{I_1}_2 
		&\leq \frac{\omegaa\big(x,\mathsf{E}\varrho_Y\big)}{\int \!\omegaa\big(x,\mathsf{E}\varrho_Y\big)\,d\varrho_X(x)} \N{x-\widehat{x}}_2
		\leq \frac{e^{-\alpha\underbar\CE(\mathsf{E}\varrho_Y)}}{\int \!\omegaa\big(x,\mathsf{E}\varrho_Y\big)\,d\varrho_X(x)} \N{x-\widehat{x}}_2
		\leq C_K^\alpha \N{x-\widehat{x}}_2.
	\end{align*}
	For $I_2$ and $I_3$, on the other hand, let us first notice that \revised{ for some $s,s'\in[0,1]$} it holds
	\begin{align*}
		&\abs{\omegaa\big(x,\mathsf{E}\varrho_Y\big)-\omegaa\big(\widehat{x},\mathsf{E}\widehat\varrho_Y\big)}
		\leq \abs{\omegaa\big(x,\mathsf{E}\varrho_Y\big)-\omegaa\big(\widehat{x},\mathsf{E}\varrho_Y\big)\big)} + \abs{\omegaa\big(\widehat{x},\mathsf{E}\varrho_Y\big)-\omegaa\big(\widehat{x},\mathsf{E}\widehat\varrho_Y\big)} \\
		&\quad\,\revised{=\abs{\partial_x\omegaa(x+s(\widehat{x}-x),\mathsf{E}\varrho_Y)}\N{x-\widehat{x}}_2} 
		\revised{+ \abs{\partial_y\omegaa(\hat x,\mathsf{E}\varrho_Y+s'(\mathsf{E}\widehat\varrho_Y-\mathsf{E}\varrho_Y))}\N{\mathsf{E}\widehat\varrho_Y-\mathsf{E}\varrho_Y}_2}\\
		&\quad\,\leq \alpha C_1e^{-\alpha\underbar\CE(\mathsf{E}\varrho_Y)}2\left(1+\N{x}_2+\Nnormal{\widehat{x}}_2+\,\Nbig{\mathsf{E}\varrho_Y}_2\right)\N{x-\widehat{x}}_2 \\
		&\quad\,\quad\,
		+ \alpha C_1e^{-\alpha\underbar\CE\left(\mathsf{E}\varrho_Y+s'(\mathsf{E}\widehat\varrho_Y-\mathsf{E}\varrho_Y)\right)}  2\left(1+\N{\mathsf{E}\widehat\varrho_Y}_2+\Nnormal{\widehat{x}}_2+\Nbig{\mathsf{E}\varrho_Y}_2\right)\N{\mathsf{E}\widehat\varrho_Y-\mathsf{E}\varrho_Y}_2
	\end{align*}
	due to Assumptions~\ref{ass:sec:wp_2} and \ref{ass:sec:wp_3}.
	With this we immediately obtain the upper bounds
	\begin{align*}
		\N{I_2}_2
		&\leq 2\alpha C_1C_K^\alpha \N{\widehat{x}}_2\left(1+\N{x}_2+\N{\mathsf{E}\widehat\varrho_Y}_2+2\Nnormal{\widehat{x}}_2+2\,\Nbig{\mathsf{E}\varrho_Y}_2\right)\left(\N{x-\widehat{x}}_2+\N{\mathsf{E}\widehat\varrho_Y-\mathsf{E}\varrho_Y}_2\right), \\
		\N{I_3}_2
		&\leq 2\alpha C_1(C_K^\alpha)^2 \N{\widehat{x}}_2\! 
		\\
		&\quad\quad  \cdot \iint \!\left(1+\N{x}_2\!+\N{\mathsf{E}\widehat\varrho_Y}_2+\!2\Nnormal{\widehat{x}}_2\!+\!2\,\Nbig{\mathsf{E}\varrho_Y}_2\right)
 \left(\N{x\!-\!\widehat{x}}_2 \!+\! \N{\mathsf{E}\widehat\varrho_Y\!-\!\mathsf{E}\varrho_Y}_2\right)d\pi(x,y,\widehat{x},\widehat{y}).
	\end{align*}
	Collecting the latter three estimates for $\N{I_1}_2,\N{I_2}_2$ and $\N{I_3}_2$ eventually gives after an application of Jensen's and Cauchy-Schwarz inequality
	\begin{align*}
		&\N{\conspointx{\varrho_X}-\conspointx{\widehat\varrho_X}}_2
		\leq C(\alpha,C_1,C_K^\alpha,K) \sqrt{\iint \N{x-\widehat{x}}_2^2 d\pi(x,y,\widehat{x},\widehat{y}) + \N{\mathsf{E}\widehat\varrho_Y-\mathsf{E}\varrho_Y}_2^2} \\
		&\qquad\leq C(\alpha,C_1,C_K^\alpha,K) \sqrt{\iint \N{x-\widehat{x}}_2^2 + \N{y-\widehat{y}}_2^2 d\pi(x,y,\widehat{x},\widehat{y})},
	\end{align*}
	where the last step is due to Jensen's inequality.
	$\Nbig{\conspointy{\varrho_Y}-\conspointy{\widehat\varrho_Y}}_2$ can be bounded analogously.
	Eventually, optimizing over all couplings $\pi\in\Pi(\varrho,\widehat\varrho)$ gives the claim.
\end{proof}


\begin{proof}[Proof of Theorem~\ref{theorem:well-posedness_mean-field}]
	The proof is based on the Leray-Schauder fixed point theorem and follows in the spirit of~\cite[Theorems~3.1, 3.2]{carrillo2018analytical}.

	\noindent\textbf{Step 1:} For given functions $u^X\in\mathcal{C}([0,T],\bbR^{d_1})$, $u^Y\in\mathcal{C}([0,T],\bbR^{d_2})$ and measure~$\rho_0\in\CP_4(\bbR^{d_1+d_2})$, according to standard SDE theory~\cite[Chapter~6]{arnold1974stochasticdifferentialequations}, we can uniquely solve the SDE system	
		\begin{subequations}\label{eq:saddle_point_dynamics_macro_tilde}
		\begin{align} 
			d\widetilde{X}_t
				&= -\lambda_1\big(\widetilde{X}_t - u_t^X\big)\,dt
				+ \sigma_1 D\big(\widetilde{X}_t - u_t^X\big)\, dB_t^{X}, \\
			d\widetilde{Y}_t
				&= -\lambda_2\big(\widetilde{Y}_t - u_t^Y\big)\,dt
				+ \sigma_2 D\big(\widetilde{Y}_t - u_t^Y\big)\,  dB_t^{Y}
		\end{align}
		\end{subequations}
		with $(\widetilde{X}_0,\widetilde{Y}_0)\sim\rho_0$ as a consequence of the coefficients being locally Lipschitz continuous and having at most linear growth.
		This induces $\widetilde\rho_t=\Law((\widetilde{X}_t, \widetilde{Y}_t))$.
		Moreover, the regularity of the initial distribution $\rho_0 \in \CP_4(\bbR^{d_1+d_2})$ allows for a fourth-order moment estimate of the form~$\mathbb{E}\big[\Nnormal{\widetilde{X}_t}_2^4+\Nnormal{\widetilde{Y}_t}_2^4\big] \leq \big(1+2\mathbb{E}\big[\Nnormal{\widetilde{X}_0}_2^4+\Nnormal{\widetilde{Y}_0}_2^4\big]\big)e^{ct}$, see, e.g.\@~\cite[Chapter~7]{arnold1974stochasticdifferentialequations}.
		So, in particular, $\widetilde\rho\in\CC([0,T],\CP_4(\bbR^{d_1+d_2}))$, i.e., $\sup_{t\in[0,T]} \iint \N{x}_2^4+\N{y}_2^4 d\widetilde\rho_t(x,y) \leq K$ for some $K>0$.
		
		\noindent\textbf{Step 2:} Let us now define the test function space
		\begin{equation} \label{eq:function_space}
		\begin{split}
			\CC^2_{*}(\bbR^{d_1+d_2})
			:=&\,\bigg\{\phi\in\CC^2(\bbR^{d_1+d_2}): 
			\,\Nnormal{\nabla\phi(x,y)}_2\leq \revised{C_\phi\big(1+\Nnormal{x}_2+\Nnormal{y}_2\big) \mbox{ for some }C_\phi>0} \\
			&\qquad\text{ and } \max\left\{\max_{k=1,\dots,d_1}\Nbig{\partial^2_{x_kx_k} \phi}_\infty, \max_{k=1,\dots,d_2}\Nbig{\partial^2_{y_ky_k} \phi}_\infty\right\} < \infty\bigg\}.
		\end{split}
		\end{equation}
		For any $\phi\in\CC^2_{*}(\bbR^d)$, by It\^o's formula, we can derive
		\begin{equation*}
		\begin{split}
			d\phi(\widetilde{X}_t,\widetilde{Y}_t)
			= &\,\nabla_x\phi(\widetilde{X}_t,\widetilde{Y}_t)\cdot\Big(-\lambda_1\big(\widetilde{X}_t - u_t^X\big)\,dt
				+ \sigma_1 D\big(\widetilde{X}_t - u_t^X\big)\, dB_t^{X}\Big)\\
			& + \nabla_y\phi(\widetilde{X}_t,\widetilde{Y}_t)\cdot\Big(-\lambda_2\big(\widetilde{Y}_t - u_t^Y\big)\,dt
				+ \sigma_2 D\big(\widetilde{Y}_t - u_t^Y\big)\,  dB_t^{Y}\Big)\\
			& + \frac{\sigma_1^2}{2}\sum_{k=1}^{d_1} \partial^2_{x_kx_k}\phi(\widetilde{X}_t,\widetilde{Y}_t) \big(\widetilde{X}_t - u_t^X\big)_{k}^2\,dt
				+ \frac{\sigma_2^2}{2}\sum_{k=1}^{d_2} \partial^2_{y_ky_k}\phi(\widetilde{X}_t,\widetilde{Y}_t) \big(\widetilde{Y}_t - u_t^Y\big)_{k}^2\,dt.
		\end{split}
		\end{equation*}
		After taking the expectation, applying Fubini's theorem and {observing that the stochastic integrals~$\bbE\int_0^t\nabla_x\phi(\widetilde{X}_t,\widetilde{Y}_t)\cdot D\big(\widetilde{X}_t - u_t^X\big)\,dB_t^{X}$ and $\bbE\int_0^t\nabla_y\phi(\widetilde{X}_t,\widetilde{Y}_t)\cdot D\big(\widetilde{Y}_t - u_t^Y\big)\,dB_t^{Y}$ vanish as a consequence of \cite[Theorem 3.2.1(iii)]{oksendal2013stochastic} due to the established regularity $\widetilde\rho\in\CC([0,T],\CP_4(\bbR^{d_1+d_2}))$ and $\phi\in\CC^2_{*}(\bbR^{d_1+d_2})$, we obtain}
		\begin{equation*}
		\begin{split}
			\frac{d}{dt}\bbE\phi(\widetilde{X}_t,\widetilde{Y}_t)
			= &\, - \lambda_1 \bbE\nabla_x\phi(\widetilde{X}_t,\widetilde{Y}_t)\cdot\big(\widetilde{X}_t - u_t^X\big)  - \lambda_2 \bbE\nabla_y\phi(\widetilde{X}_t,\widetilde{Y}_t)\cdot\big(\widetilde{Y}_t - u_t^Y\big) \\
			&\, + \frac{\sigma_1^2}{2} \bbE\sum_{k=1}^{d_1} \partial^2_{x_kx_k}\phi(\widetilde{X}_t,\widetilde{Y}_t) \big(\widetilde{X}_t - u_t^X\big)_{k}^2
				+ \frac{\sigma_2^2}{2} \bbE\sum_{k=1}^{d_2} \partial^2_{y_ky_k}\phi(\widetilde{X}_t,\widetilde{Y}_t) \big(\widetilde{Y}_t - u_t^Y\big)_{k}^2
		\end{split}
		\end{equation*}
		as a consequence of the fundamental theorem of calculus.
		This shows that the measure $\widetilde\rho\in\CC([0,T],\CP_4(\bbR^{d_1+d_2}))$ satisfies the Fokker-Planck equation
		\begin{align} \label{proof:auxiliary_Fokker-Planck}
		\begin{aligned}
			&\frac{d}{dt}\int\!\phi(x,y)\,d\widetilde\rho_t(x,y) =
			 - \!\int \!\lambda_1 \nabla_x\phi(x,y)\cdot\big(x \!-\! u_t^X\big) \!+\! \lambda_2 \nabla_y\phi(x,y)\cdot\big(y \!-\! u_t^Y\big)\,d\widetilde\rho_t(x,y) \\
			&\qquad\quad +\! \int \!\frac{\sigma_1^2}{2} \sum_{k=1}^{d_1} \partial^2_{x_kx_k}\phi(x,y) \big(x \!-\! u_t^X\big)_{k}^2
				\!+\! \frac{\sigma_2^2}{2} \sum_{k=1}^{d_2} \partial^2_{y_ky_k}\phi(x,y) \big(y \!-\! u_t^Y\big)_{k}^2 \,d\widetilde\rho_t(x,y).
		\end{aligned}
		\end{align}
		\noindent\textbf{Step 3:} Setting $\CT u:=\big((\conspointx{\widetilde\rho_X})^T,(\conspointy{\widetilde\rho_Y})^T\big)^T\in\mathcal{C}([0,T],\bbR^{d_1+d_2})$ provides the self-mapping property of the map
		\begin{align*}
			\CT:\CC([0,T],\bbR^{d_1+d_2})\rightarrow\CC([0,T],\bbR^{d_1+d_2}), \quad u\mapsto\mathcal{T}u=\big((\conspointx{\widetilde\rho_X})^T,(\conspointy{\widetilde\rho_Y})^T\big)^T.
		\end{align*}
		By means of Lemma~\ref{lem:stabiliy_conspoint} we have
		\begin{align*}
			\N{\conspointx{\widetilde\rho_{X,t}}-\conspointx{\widetilde\rho_{X,s}}}_2 + \N{\conspointy{\widetilde\rho_{Y,t}}-\conspointy{\widetilde\rho_{Y,s}}}_2
		\leq c_0 W_2(\widetilde\rho_{t},\widetilde\rho_{s})
		\lesssim c_0 \abs{t-s}^{\frac{1}{2}},
		\end{align*}
		which shows the H\"older-$1/2$ continuity of $\CT$ due to the compact embedding $\CC^{0,1/2}([0,T],\bbR^{d_1+d_2}) \hookrightarrow \CC([0,T],\bbR^{d_1+d_2})$.
		The last inequality, we note that due to It\^o's isometry it holds
		\begin{align*}
			W_2^2(\widetilde\rho_{t},\widetilde\rho_{s})
			&\!\leq\! \bbE\!\left[\Nbig{\widetilde{X}_t\!-\!\widetilde{X}_s}_2^2 \!+\! \Nbig{\widetilde{Y}_t\!-\!\widetilde{Y}_s}_2^2 \right] 
			\!\leq\! 4\!\left((\lambda_1^2\!+\!\lambda_2^2)T\!+\!(\sigma_1^2\!+\!\sigma_2^2)\right)\! \left(K \!+\! \Nbig{u^X}_{L^\infty}^2 \!\!+\!\Nbig{u^Y}_{L^\infty}^2\!\right)\!\abs{t\!-\!s}\!.
		\end{align*}
		
		\noindent\textbf{Step 4:} Now, for $u\in\mathcal{C}([0,T],\bbR^{d_1+d_2})$ satisfying $u=\vartheta\CT u$ with $\vartheta\in[0,1]$, there exists $\rho\in\CC([0,T],\CP_4(\bbR^{d_1+d_2}))$ satisfying~\eqref{proof:auxiliary_Fokker-Planck} such that $u=\vartheta \big((\conspointx{\rho_X})^T,(\conspointy{\rho_Y})^T\big)^T$.
		As a consequence of Lemma~\ref{lem:aux_lemma_bound} with $s=0$ we can show that
		\begin{align*}
			\N{u_t}_2^2
			&\leq \vartheta^2 \left((C_K^\alpha)^2\int \N{x}_2^2 d\rho_{X,t}(x) + (C_K^\beta)^2 \int \N{y}_2^2 d\rho_{Y,t}(y)\right)
			\leq \vartheta^2 \left((C_K^\alpha)^2 + (C_K^\beta)^2\right) \sqrt{K}
		\end{align*}
		where we can use $K=\big(1+2\mathbb{E}\big[\Nnormal{\widetilde{X}_0}_2^4+\Nnormal{\widetilde{Y}_0}_2^4\big]\big)e^{cT}$ as of Step~$1$.
		This allows for a uniform estimate of $\N{u}_{L^\infty}<q$ for $q>0$.
		An application of the Leray-Schauder fixed point theorem concludes the proof of existence by providing a solution to~\eqref{eq:saddle_point_dynamics_macro}.
		
		\noindent\textbf{Step 5:} For uniqueness, suppose we have two fixed points~$u^1$ and $u^2$ (as specified in the previous step) together with corresponding processes~$((\widetilde{X}^{1})^T,(\widetilde{Y}^{1})^T)^T$ and $((\widetilde{X}^{2})^T,(\widetilde{Y}^{2})^T)^T$ satisfying~\eqref{eq:saddle_point_dynamics_macro_tilde}.
		Then, taking their difference while keeping the initial conditions and respective Brownian motion paths, we obtain after an application of It\^o's isometry and employment of Lemma~\ref{lem:stabiliy_conspoint} the bound
		\begin{align*}
			&\bbE\left[\Nbig{\widetilde{X}^1_t-\widetilde{X}^2_t}_2^2 + \Nbig{\widetilde{Y}^2_t-\widetilde{Y}^2_t}_2^2 \right]\\
			&\qquad \leq c\bbE\int_0^t \Nbig{\widetilde{X}^1_\tau\!-\!\widetilde{X}^2_\tau}_2^2 \!+\! \Nbig{\widetilde{Y}^1_\tau\!-\!\widetilde{Y}^2_\tau}_2^2 \!+\! \N{\conspointx{\widetilde\rho^1_{X,\tau}}\!-\!\conspointx{\widetilde\rho^2_{X,\tau}}}_2^2 \!+\! \N{\conspointy{\widetilde\rho^1_{Y,\tau}}\!-\!\conspointy{\widetilde\rho^2_{Y,\tau}}}_2^2\, d\tau \\
			&\qquad \lesssim c \bbE\int_0^t \Nbig{\widetilde{X}^1_\tau\!-\!\widetilde{X}^2_\tau}_2^2 \!+\! \Nbig{\widetilde{Y}^1_\tau\!-\!\widetilde{Y}^2_\tau}_2^2 \, d\tau
		\end{align*}
		with $c=4\left((\lambda_1^2+\lambda_2^2)T+(\sigma_1^2+\sigma_2^2)\right)$.
		Gr\"onwall's inequality eventually shows uniqueness since $\bbE\big[\Nnormal{\widetilde{X}^1_t-\widetilde{X}^2_t}_2^2 + \Nnormal{\widetilde{Y}^2_t-\widetilde{Y}^2_t}_2^2 \big]=0$ for all $t\in[0,T]$.
\end{proof}

\section{Convergence to Saddle Points} \label{sec:convergence}
Inspired by the theories of mean-field limits for consensus-based optimization methods (see \cite{fornasier2021consensus,huang2021MFLCBO,fornasier2020consensus_hypersurface_wellposedness} for instance), the convergence of particles systems \eqref{eq:saddle_point_dynamics_micro} to the mean-field dynamics~\eqref{eq:saddle_point_dynamics_macro} follows in a similar way and thus the associated argument is omitted. In this section we present the main theoretical result of our paper concerned with the convergence of the macroscopic dynamics~\eqref{eq:saddle_point_dynamics_macro} towards saddle points of objective functions~$\CE$ that satisfy the following conditions.

\begin{definition}[Assumptions] \label{ass:sec:convergence}
	Throughout this section we are interested in objective functions $\CE\in\CC^2(\bbR^{d_1+d_2})$, for which
	\begin{enumerate}[label=A\arabic*,labelsep=10pt,leftmargin=35pt]
		\item \label{ass:sec:convergence_1}
			there exist two functions~$\underbar{\CE}\in\CC^1(\bbR^{d_2})$ and $\overbar{\CE}\in\CC^1(\bbR^{d_1})$ such that
			\begin{equation*}
				\underbar{\CE}(y)\leq\CE(x,y)\leq\overbar{\CE}(x)
			\end{equation*}
			for all $(x,y)\in\bbR^{d_1+d_2}$.
		The functions~$\underbar{\CE}$ and~$\overbar{\CE}$ shall, for some constant~$\overbar{C}_{\nabla\CE}>0$, satisfy $\Nnormal{\nabla\underbar{\CE}(y)}_2\leq \overbar{C}_{\nabla\CE}$ for all $y\in\bbR^{d_2}$ and $\Nnormal{\nabla\overbar{\CE}(x)}_2\leq \overbar{C}_{\nabla\CE}$ for all $x\in\bbR^{d_1}$.
		\item \label{ass:sec:convergence_2}
			there exist constants~$C_{\nabla\CE},C_{\nabla^2\CE}>0$ such that  
			\begin{equation*}
			\begin{split}
				&\max\left\{
					\sup_{(x,y)\in\bbR^{d_1}\times \bbR^{d_2}}\N{\nabla_x\CE(x,y)}_2,
					\sup_{(x,y)\in\bbR^{d_1}\times \bbR^{d_2}}\N{\nabla_y\CE(x,y)}_2
				\right\}
				\leq C_{\nabla\CE},
				\text{ and} \\
				&\max\left\{
					\max_{k=1,\dots,d_1} \N{\partial^2_{x_kx_k}\CE}_\infty,
					\max_{k=1,\dots,d_2} \N{\partial^2_{y_ky_k}\CE}_\infty,
					\N{\rho\left(\nabla_x^2\CE\right)}_\infty,
					\N{\rho\left(\nabla_y^2\CE\right)}_\infty
				\right\}
				\leq C_{\nabla^2\CE},
			\end{split}
			\end{equation*}
			where $\N{\,\cdot\,}_\infty$ denotes the $L^\infty$ norm on $\CC\left(\bbR^{d_1+d_2}\right)$ and $\rho$ denotes the spectral radius.
		\item \label{ass:sec:ICP}
			there exist constants~$\epsilon_0,\eta,\nu>0$ such that \revised{for each $(x,y)\in\bbR^{d_1+d_2}$ satisfying 
		$	\CE^*-\CE(x^*, {y}) \leq \epsilon_0$
		  and 
		$\CE( {x},y^*)-\CE^* \leq \epsilon_0 $ for some
		 saddle point $\left(\saddlepoint\right)$ of $\CE$,} we have			\begin{equation*}
				\N{x-x^*}_2
				\leq \frac{1}{\eta}\big(\absnormal{\CE(x,y^*)-\CE^*}\big)^\nu \mbox{ and }\N{y-y^*}_2
				\leq \frac{1}{\eta}\big(\absnormal{\CE(x^*,y)-\CE^*)}\big)^\nu.
			\end{equation*}
	\end{enumerate}
\end{definition}

Assumption~\ref{ass:sec:convergence_1} requires that the twice continuously differentiable objective function~$\CE$ is bounded from below by a function $\underbar{\CE}$, which depends only on the $y$ coordinate, and from above by a function $\overbar{\CE}$, which depends only on the $x$ coordinate.
Moreover, the first order derivatives of $\underbar{\CE}$ and $\overbar{\CE}$ are assumed to be uniformly bounded.

Assumption~\ref{ass:sec:convergence_2} is a mere technical regularity assumption about $\CE$ in terms of the first and second derivatives.
In particular, it requires that the gradients as well as second order derivatives of $\CE$ are uniformly bounded, which is however necessary only for theoretical analysis of the long time behavior of the algorithm.
Analogous regularity assumptions may be found in the literature, see, e.g., \cite{carrillo2019consensus,carrillo2018analytical,qiu2022PSOconvergence}.
However, as a purely zero-order derivative-free method, our CBO-SP algorithm only uses point-wise values of the objective function~$\CE$ in practical applications.

Assumption~\ref{ass:sec:ICP} on the other hand should be regarded as a tractability condition on the landscape of the objective function~$\CE$.
It imposes coercivity of $\CE$ around saddle points, which relates the distance from $(\saddlepoint)$ with the value of the objective function.
We refer to the discussion after \cite[Remark~9]{fornasier2021consensus} for related conditions in the machine learning literature.

\revised{In order to formulate in Theorem~\ref{theorem:convergece} below the result about the convergence of the dynamics~\eqref{eq:saddle_point_dynamics_macro} towards saddle points of the objective functions~$\CE$ satisfying the aforementioned assumption, let us define the variances
\begin{equation} \label{eq:variances}
	\VarX(t) = \bbE\N{\OX_t-\bbE\OX_t}_2^2
	\quad\text{and}\quad
	\VarY(t) = \bbE\N{\OY_t-\bbE\OY_t}_2^2,
\end{equation}
which act as Lyapunov functionals of the dynamics.}
In addition, we require certain well-preparedness assumptions about the initial data and parameters.
For this reason we introduce the notations 
\begin{subequations} \label{eq:CMs}
\begin{align}
	\widetilde\CM^X(t) &:= \bbE\exp\left(-\alpha\CE(\OX_t,\bbE\OY_t)\right)
	&&\text{and} \qquad \
	\widetilde\CM^Y(t) := \bbE\exp\left(\beta\CE(\bbE\OX_t,\OY_t)\right) \label{eq:CM_tilde},\\
	\CM^X(t) &:= \widetilde\CM^X(t)\,e^{\alpha\underbar{\CE}(\bbE\overbarscript{Y}_t)}
	&&\text{and} \qquad \
	\CM^Y(t) := \widetilde\CM^Y(t)\,e^{-\beta\overbar{\CE}(\bbE\overbarscript{X}_t)} \label{eq:CM},\\
	\CM^X_*(t) &:= \bbE\exp\left(-\alpha\CE(\OX_t,y^*)\right)
	&&\text{and} \qquad \
	\CM^Y_*(t) := \bbE\exp\left(\beta\CE(x^*,\OY_t)\right), \label{eq:CM_star}
\end{align}
\end{subequations}
\revised{where the definitions of $\CM^X_*$ and $\CM^Y_*$ in \eqref{eq:CM_star} may be different for potentially different saddle points $(x^*,y^*)$ (since our assumptions allow for the existence of multiple such points), meaning, in particular, that there may be multiple functionals $\CM^X_*$ and $\CM^Y_*$.
This ambiguity is clarified in Remark~\ref{remark_multiple_saddlepoints}.}

\begin{definition}[Well-preparedness of the initial data and parameters] \label{ass:sec:convergence}	
	The initial datum $(\OX_0,\OY_0)$ and the parameters~$\alpha,\beta,\lambda_1,\lambda_2,\sigma_1$ and $\sigma_2$ of the CBO-SP method are well-prepared if  
	\begin{enumerate}[label=P\arabic*,labelsep=10pt,leftmargin=35pt]
		\item \label{asm:well_preparedness_P1}
		$\mu_1:= 2\left(\lambda_1 - 4\sigma_1^2/\CM^X(0)\right)>0$ and $\mu_2:= 2\left(\lambda_2 - 4\sigma_2^2/\CM^Y(0)\right)>0$,
\item \label{asm:well_preparedness_P1'}
\revised{all saddle points $(x^*,y^*)$ lie in $\supp(\rho_0)$, where $\rho_0:=\Law(\OX_0,\OY_0)$ has marginals $\rho_0^X$ and $\rho_0^Y$.	
	Moreover, for any $\delta>0$, there exists some constant $C_{\delta}>0$ depending only on $\delta$ such that it holds
	\begin{align*}
	\begin{aligned}
		\rho_0^X\left(\left\{x:\exp\left(-\CE(x,\bbE(\OY_0))\right)>\exp\left(-\min_{x\in\bbR^{d_1}}\CE(x,\bbE(\OY_0))\right)-\delta \right\}\right)
		&\geq C_\delta, \\
		\rho_0^Y\left(\left\{y:\exp\left(\CE(\bbE(\OX_0),y)\right)>\exp\left(\max_{y\in\bbR^{d_2}}\CE(\bbE(\OX_0),y)\right)-\delta \right\}\right)
		&\geq C_\delta
	\end{aligned}
	\end{align*}
	as well as
	\begin{align*}
	\begin{aligned}
		\rho_0^X\big(\!\left\{x:\exp\left(-\CE(x,y^*)\right)>\exp\left(-\CE^*\right)-\delta \right\}\!\big)
		&\geq C_\delta, \\
		\rho_0^Y\big(\!\left\{y:\exp\left(\CE(x^*,y)\right)>\exp\left(\CE^*\right)-\delta \right\}\!\big)
		&\geq C_\delta.
	\end{aligned}
	\end{align*}
}
		\item it \revised{holds that }\label{asm:well_preparedness_P2}
		\begin{align*} 
			&4\alpha C_{\nabla^2\CE} \left(\lambda_1 + \frac{\sigma_1^2}{2}\right) \frac{\VarX(0)}{\mu_1\CM^X(0)}
			+2\sqrt{2}\alpha\lambda_2 C_{\nabla\CE} \frac{\sqrt{\VarY(0)}}{\mu_2\sqrt{\CM^Y(0)}}
			\leq\frac{1}{8}\CM^X(0),\\
			&4\beta C_{\nabla^2\CE} \left(\lambda_2 + \frac{\sigma_2^2}{2}\right) \frac{\VarY(0)}{\mu_2\CM^Y(0)}
			+2\sqrt{2}\beta\lambda_1 C_{\nabla\CE} \frac{\sqrt{\VarX(0)}}{\mu_1\sqrt{\CM^X(0)}}
			\leq\frac{1}{8}\CM^Y(0),
		\end{align*}
		\item it \revised{holds for any fixed $(x^*,y^*)$ that} \label{asm:well_preparedness_P3}
		\begin{align*} 
			&4\alpha \sigma_1^2C_{\nabla^2\CE} \left(\lambda_1 + \frac{\sigma_1^2}{2}\right)   \frac{\VarX(0)}{\mu_1\CM^X(0)}
			+8\alpha\lambda_1 C_{\nabla\CE}\frac{\sqrt{\VarX(0)}}{\mu_1\sqrt{\CM^X(0)}}
			+\sqrt{2}\alpha\lambda_2 C_{\nabla\CE} \frac{\sqrt{\VarY(0)}}{\mu_2\sqrt{\CM^Y(0)}}\\
			&\qquad\qquad\qquad\quad
			\leq\frac{1}{8}\min\big\{\CM^X_*(0)e^{\alpha\underbar{\CE}(y^*)},\widetilde\CM^X(0)e^{-\alpha\CE_M}\big\},\\
			&4\beta C_{\nabla^2\CE} \left(\lambda_2 + \frac{\sigma_2^2}{2}\right) \frac{\VarY(0)}{\mu_2\CM^Y(0)}
			+8\beta\lambda_2 C_{\nabla\CE}\frac{\sqrt{\VarY(0)}}{\mu_2\sqrt{\CM^Y(0)}}
			+\sqrt{2}\beta\lambda_1 C_{\nabla\CE} \frac{\sqrt{\VarX(0)}}{\mu_1\sqrt{\CM^X(0)}}\\
			&\qquad\qquad\qquad\quad
			\leq\frac{1}{8}\min\big\{\CM^Y_*(0)e^{-\beta\overbar{\CE}(x^*)},\widetilde\CM^Y(0)e^{-\beta\CE_M}\big\},
		\end{align*}
		where $\CE_M$ depends on $(\widetilde{x},\widetilde{y})$ as well as $\underbar{\CE}$ and $\overbar{\CE}$ (see Theorem \ref{theorem:convergece} below).
	\end{enumerate}
\end{definition}

Condition~\ref{asm:well_preparedness_P1} can be satisfied appropriately big choices of $\lambda_1$ and $\lambda_2$.
\revised{Condition~\ref{asm:well_preparedness_P1'} is valid if the initial distribution~$\rho_0$ has some mass at the saddle points $(x^*,y^*)$.
While this may have a certain locality flavor, in the case of functions $\CE$ having multiple saddle points, the condition is generically satisfied at least for one of them, allowing essentially to obtain a global result.
It is actually sufficient if there is at least one saddle point satisfying Condition~\ref{asm:well_preparedness_P1'}, in which case the method is agnostic to saddle points not in $\supp(\rho_0)$.}
Conditions~\ref{asm:well_preparedness_P2}--\ref{asm:well_preparedness_P3} on the other hand may be ensured if the initial variances~$\VarX(0)$ and $\VarY(0)$ are sufficiently small.
Well-preparedness conditions similar to \ref{asm:well_preparedness_P1'}--\ref{asm:well_preparedness_P3} can be found in the literature about a convergence analysis of CBO for minimization, see, e.g., \cite{carrillo2018analytical,carrillo2019consensus,fornasier2020consensus_sphere_convergence,qiu2022PSOconvergence}, while we note that the coupling of $(\OX,\OY)$ due to the intrinsic difference between games and optimizations prompts us to use some different techniques in the proof.

\revised{We are now ready to state the result about the convergence of the dynamics~\eqref{eq:saddle_point_dynamics_macro} towards saddle points of the objective functions~$\CE$.}
The proof details are deferred to Section~\ref{sec:proofs}.


\revised{\begin{theorem} \label{theorem:convergece}
	Let $\CE$ satisfy Assumptions~\ref{ass:sec:convergence_1} and \ref{ass:sec:convergence_2} and let $(\OX_t,\OY_t)_{t\geq0}$ be a solution to the SDE~\eqref{eq:saddle_point_dynamics_macro}. Then the following statements hold.
	\begin{enumerate}[label=(\arabic*),labelsep=10pt,leftmargin=35pt]
		\item \label{theorem:convergece_item1} Under the assumption of well-preparedness of the initial datum $(\OX_0,\OY_0)$ and the parameters~$\alpha,\beta,\lambda_1,\lambda_2,\sigma_1$ and $\sigma_2$ in the sense of \ref{asm:well_preparedness_P1}--\ref{asm:well_preparedness_P2},
		$\VarX$ and $\VarY$ as defined in~\eqref{eq:variances} converge exponentially fast to $0$ as $t\rightarrow\infty$.
		More precisely, it holds
		\begin{equation} \label{eq:theorem:convergece:3}
			\VarX(t)+\VarY(t) \leq \VarX(0)e^{-\mu_1t}+\VarY(0)e^{-\mu_2t}.
		\end{equation}
		Moreover, there exists some $(\widetilde{x},\widetilde{y})$ depending in particular on $\alpha$ and $\beta$ such that 
		\begin{equation} \label{eq:theorem:convergece:4}
			\left(\bbE\OX_t,\bbE\OY_t) \rightarrow (\widetilde{x},\widetilde{y}\right) 
			\quad\text{and}\quad
			\left(\conspointx{\rho_{X,t}},\conspointy{\rho_{Y,t}}\right) \rightarrow (\widetilde{x},\widetilde{y})
		\end{equation}
		as $t\rightarrow\infty$. 
	\item \label{theorem:convergece_item2} For any given accuracy $\varepsilon>0$, there exist some $\alpha_0,\beta_0>0$ such that for all $\alpha\geq\alpha_0$ and $\beta\geq\beta_0$ the point $(\widetilde{x},\widetilde{y})$ from \ref{theorem:convergece_item1} (which may depend on $\alpha$ and $\beta$) satisfies
	\begin{equation} \label{eq:theorem:convergece:5}
		\absnormal{\CE(\widetilde{x},\widetilde{y})-\CE^*} \leq \varepsilon
		\quad\text{as well as}\quad
		\CE^*-\CE(x^*,\widetilde{y}) \leq \varepsilon
		\quad\text{and}\quad
		\CE(\widetilde{x},y^*)-\CE^* \leq \varepsilon
	\end{equation}
	provided that the well-preparedness Assumptions~\ref{asm:well_preparedness_P1}--\ref{asm:well_preparedness_P3} hold for such $\alpha$ and $\beta$ together with the initial datum~$(\OX_0,\OY_0)$.
\item \label{theorem:convergece_item3} If $\CE$ satisfies Assumption~\ref{ass:sec:ICP} with respect to $(\widetilde{x},\widetilde{y})$ from \ref{theorem:convergece_item2} with $\varepsilon\leq \epsilon_0$, i.e., there exists some saddle point $(x^*,y^*)$ depending on $(\widetilde{x},\widetilde{y})$ such that
	\begin{equation*}
		\|\widetilde{x}-x^*\|_2\leq \frac{1}{\eta}\left(\abs{\mathcal{E}(\widetilde{x},y^*)-\mathcal{E}^*}\right)^\nu\quad \mbox{ and }\quad \|\widetilde{y}-y^*\|_2\leq \frac{1}{\eta}\left(\abs{\mathcal{E}(x^*,\widetilde{y})-\mathcal{E}^*}\right)^\nu,
	\end{equation*}
	then we have
	\begin{equation} \label{eq:theorem:convergece:7}
		\N{\left(\widetilde{x},\widetilde{y}\right)-\left(\saddlepoint\right)}_2
		\leq \frac{2}{\eta} \varepsilon^\nu
	\end{equation}
	provided that the well-preparedness Assumptions~\ref{asm:well_preparedness_P1}--\ref{asm:well_preparedness_P3} hold for sufficiently large $\alpha$ and $\beta$ together with the initial datum~$(\OX_0,\OY_0)$.
	\end{enumerate}
\end{theorem}}

The \revised{first part~\ref{theorem:convergece_item1}} of Theorem~\ref{theorem:convergece} states that under suitable well-preparedness conditions on the initialization and the parameters, the mean-field dynamics~\eqref{eq:saddle_point_dynamics_macro} reaches consensus at \revised{\textit{some}} location $(\widetilde{x},\widetilde{y})$, \revised{which may depend in particular on $\alpha$ and $\beta$,} as time evolves.
In the \revised{second part~\ref{theorem:convergece_item2} of the statement, for sufficiently large $\alpha$ and $\beta$ as well as under certain well-preparedness conditions, properties of \textit{the} corresponding point~$(\widetilde{x},\widetilde{y})$ are specified which are typical for saddle points, see~\eqref{eq:theorem:convergece:5}.}
\revised{These properties eventually allow to conclude in the third part~\ref{theorem:convergece_item3} of the result, that \textit{the} $(\widetilde{x},\widetilde{y})$ from before is arbitrarily close to \textit{any} saddle point $(x^*,y^*)$ which satisfies the inverse continuity property~\ref{ass:sec:ICP}.}

\revised{
\begin{remark} \label{remark_multiple_saddlepoints}
	It is worth mentioning at this point that in order to prove \eqref{eq:theorem:convergece:5}, any saddle point $(x^*,y^*)$ satisfying Assumption~\ref{asm:well_preparedness_P3} can be used in the definitions of $\CM_*^X$ and $\CM_*^Y$.
	For the proof of \eqref{eq:theorem:convergece:7}, on the other hand, it is necessary to use in the definitions of $\CM_*^X$ and $\CM_*^Y$ a specific saddle point $(x^*,y^*)$ that satisfies the inverse continuity property~\ref{ass:sec:ICP} with respect to $(\widetilde{x}, \widetilde{y})$ as well as Assumption~\ref{asm:well_preparedness_P3}, so, in this case, the saddle point $(x^*,y^*)$ does depend on $(\widetilde{x}, \widetilde{y})$.
\end{remark}}

\section{Proof Details for Section~\ref{sec:convergence}} \label{sec:proofs}

In this section we provide the proof details for the convergence result of the mean-field dynamics~\eqref{eq:saddle_point_dynamics_macro} to a saddle point of the objective function~$\CE$.
Sections~\ref{subsec:evolution_variances}--\ref{subsec:evolution_CM*} present individual statements which are necessary in the proof of our main theorem, Theorem~\ref{theorem:convergece}, which is then given in Section~\ref{subsec:proof_main_theorem}.

\subsection{Time-Evolution of the Variances $\VarX$ and $\VarY$} \label{subsec:evolution_variances}

In order to ensure consensus formation of the mean-field dynamics~\eqref{eq:saddle_point_dynamics_macro} we show that the variances~$\VarX(t) = \bbE\Nnormal{\OX_t-\bbE\OX_t}_2^2$ and $\VarY(t) = \bbE\Nnormal{\OY_t-\bbE\OY_t}_2^2$ of the particle distribution decay to $0$ as $t\rightarrow\infty$.
For this we need to analyze their time-evolutions, which is the content of the following statement.

\begin{lemma} \label{lem:evolution_Var}
	Let $\VarX$ and $\VarY$ be as defined in~\eqref{eq:variances} and let us recall the definitions of $\CM^X$ and $\CM^Y$ from \eqref{eq:CM}.
	Then, under Assumption~\ref{ass:sec:convergence_1}, it \revised{holds}
	\begin{equation} \label{eq:evolution_Var}
	\begin{split}
		\frac{d}{dt}\VarX(t) \leq -2\bigg(\lambda_1 \!-\! \frac{2\sigma_1^2}{\CM^X(t)}\bigg) \VarX(t)
		\quad\text{and}\quad
		\frac{d}{dt}\VarY(t) \leq -2\bigg(\lambda_2 \!-\! \frac{2\sigma_2^2}{\CM^Y(t)}\bigg) \VarY(t).
	\end{split}
	\end{equation}
\end{lemma}

\begin{proof}
	By means of It\^o's calculus we have
	\begin{equation*} 
		d\N{\OX_t-\bbE\OX_t}_2^2
		= 2\left(\OX_t-\bbE\OX_t\right) \cdot d\OX_t + \sigma_1^2 \N{\OX_t - \conspointx{\rho_{X,t}}}_2^2dt.
	\end{equation*}
	Since the appearing stochastic integral vanishes when taking the expectation as a consequence of the regularity established in Theorem~\ref{theorem:well-posedness_mean-field} and Assumption~\ref{ass:sec:convergence_1}, we obtain
	\begin{equation} \label{eq:lem:evolution_Var_2}
	\begin{split}
		\frac{d}{dt}\VarX(t)
		&= -2\lambda_1\bbE\left[\left(\OX_t-\bbE\OX_t\right) \cdot \left(\OX_t - \conspointx{\rho_{X,t}}\right)\right] + \sigma_1^2 \bbE\N{\OX_t - \conspointx{\rho_{X,t}}}_2^2\\
		&= -2\lambda_1\VarX(t) + \sigma_1^2 \bbE\N{\OX_t - \conspointx{\rho_{X,t}}}_2^2,
	\end{split}
	\end{equation}
	where we used that $\bbE\left[\left(\OX_t-\bbE\OX_t\right) \cdot \left(\bbE\OX_t - \conspointx{\rho_{X,t}}\right)\right]=0$.
	Analogously, we derive
	\begin{equation} \label{eq:lem:evolution_Var_3}
	\begin{split}
		\frac{d}{dt}\VarY(t)
		&= -2\lambda_2\VarY(t) + \sigma_2^2 \bbE\N{\OY_t - \conspointy{\rho_{Y,t}}}_2^2.
	\end{split}
	\end{equation}
	In order to control the terms $\bbE\Nbig{\OX_t - \conspointx{\rho_{X,t}}}_2^2$ and $\bbE\Nbig{\OY_t - \conspointy{\rho_{Y,t}}}_2^2$ appearing in \eqref{eq:lem:evolution_Var_2} and \eqref{eq:lem:evolution_Var_3}, let us first observe that, for any~$\widehat{x}\in\bbR^{d_1}$, $\widehat{y}\in\bbR^{d_2}$, Jensen's inequality gives
	\begin{align}
		\N{\widehat{x}-\conspointx{\rho_{X,t}}}_2^2
		&\leq \frac{1}{\bbE\omegaa(\OX_t,\bbE\OY_t)} \int \! \N{\widehat{x}-x}_2^2 \omegaa(x,\bbE\OY_t)\,d\rho_{X,t}(x), \label{eq:lem:evolution_Var_4} \\
		\N{\widehat{y}-\conspointy{\rho_{Y,t}}}_2^2
		&\leq \frac{1}{\bbE\omegab(\bbE\OX_t,\OY_t)} \int \! \N{\widehat{y}-y}_2^2 \omegab(\bbE\OX_t,y)\,d\rho_{Y,t}(y). \label{eq:lem:evolution_Var_5}
	\end{align}
	Exploiting the boundedness of $\CE$ as of Assumption~\ref{ass:sec:convergence_1}, the two latter bounds in particular imply
	\begin{equation} \label{eq:lem:evolution_Var_6}
	\begin{split}
		\bbE\N{\OX_t\!-\!\conspointx{\rho_{X,t}}}_2^2
		&\!\leq\! \frac{2}{\bbE\omegaa(\OX_t,\bbE\OY_t)} \int \! \left(\bbE\N{\OX_t-\bbE\OX_t}_2^2 \!+\! \N{\bbE\OX_t\!-\!x}_2^2\right) \omegaa(x,\bbE\OY_t)\,d\rho_{X,t}(x) \\
		&\!\leq\! 2\VarX(t) \!+\! \frac{2}{\bbE\omegaa(\OX_t,\bbE\OY_t)}e^{-\alpha\underbar{\CE}(\bbE\overbarscript{Y}_t)} \VarX(t)
		\!\leq\! 4 \frac{\VarX(t)}{\CM^X(t)} 
	\end{split}
	\end{equation}
	and analogously
	\begin{equation} \label{eq:lem:evolution_Var_7}
		\bbE\N{\OY_t\!-\!\conspointy{\rho_{Y,t}}}_2^2
		\!\leq\! 4 \frac{\VarY(t)}{\CM^Y(t)},
	\end{equation}
	which allow to conclude the proof when being inserted into~\eqref{eq:lem:evolution_Var_2} and \eqref{eq:lem:evolution_Var_3}, respectively.
\end{proof}

\subsection{Time-Evolution of the Functionals $\CM^X$ and $\CM^Y$ \revised{from \eqref{eq:CM}}} \label{subsec:evolution_CM}

In the time-evolutions~\eqref{eq:evolution_Var} of the variances~$\VarX$ and $\VarY$ there appear the functionals~$\CM^X$ and $\CM^Y$ \revised{as defined in \eqref{eq:CM},} which need to be controlled in order to ensure that the decay rates can be bounded from below by a positive constant, which eventually leads to at least exponential decay of the variances and therefore consensus of the dynamics~\eqref{eq:saddle_point_dynamics_macro}.
We therefore investigate the evolutions of $\CM^X$ and $\CM^Y$ next.
To do so, let us recall from \eqref{eq:CM} that $\CM^X(t)=\widetilde\CM^X(t)\,e^{\alpha\underbar{\CE}(\bbE\overbarscript{Y}_t)}$ and $\CM^Y(t)=\widetilde\CM^Y(t)\,e^{-\beta\overbar{\CE}(\bbE\overbarscript{X}_t)}$.
We first bound in Lemma~\ref{lem:evolution_CM_aux} the evolutions of $\widetilde\CM^X$ and $\widetilde\CM^Y$\revised{as defined in \eqref{eq:CM_tilde},} before we use product rule to obtain a lower bound for the evolutions of $\CM^X$ and $\CM^Y$ in Lemma~\ref{lem:evolution_CM}.

Let us furthermore remark, that $\widetilde\CM^X$ and $\widetilde\CM^Y$  will later allow to characterize the convergence point of the dynamics~\eqref{eq:saddle_point_dynamics_macro}.

\begin{lemma} \label{lem:evolution_CM_aux}
	Let $\VarX$ and $\VarY$ be as defined in~\eqref{eq:variances}, and $\widetilde\CM^X$ and $\widetilde\CM^Y$ as in~\eqref{eq:CM_tilde}.
	Then, under Assumptions~\ref{ass:sec:convergence_1} and \ref{ass:sec:convergence_2}, it \revised{holds}
	\begin{equation}
	\begin{split}
		\frac{d}{dt}\widetilde\CM^X(t)
		&\geq -4\alpha e^{-\alpha\underbar{\CE}(\bbE\overbarscript{Y}_t)} C_{\nabla^2\CE}\left(\lambda_1 + \frac{\sigma_1^2}{2}\right) \frac{\VarX(t)}{\CM^X(t)} -\alpha\lambda_2 e^{-\alpha\underbar{\CE}(\bbE\overbarscript{Y}_t)} C_{\nabla\CE} \frac{\sqrt{\VarY(t)}}{\sqrt{\CM^Y(t)}}
	\end{split}
	\end{equation}
	as well as
	\begin{equation}
	\begin{split}
		\frac{d}{dt}\widetilde\CM^Y(t)
		&\geq -4\beta e^{\beta\overbar{\CE}(\bbE\overbarscript{X}_t)} C_{\nabla^2\CE} \left(\lambda_2 + \frac{\sigma_2^2}{2}\right) \frac{\VarY(t)}{\CM^Y(t)} -\beta\lambda_1 e^{\beta\overbar{\CE}(\bbE\overbarscript{X}_t)} C_{\nabla\CE} \frac{\sqrt{\VarX(t)}}{\sqrt{\CM^X(t)}}.
	\end{split}
	\end{equation}
\end{lemma}

\begin{proof}
	With It\^o's formula and chain rule we first note that
	\begin{equation*} 
	\begin{split}
		&d\widetilde\CM^X(t) = 
		-\alpha \bbE\left[\exp\left(-\alpha\CE(\OX_t,\bbE\OY_t)\right)\nabla_x\CE(\OX_t,\bbE\OY_t) \cdot d\OX_t\right]\\
		&\;\; +\frac{\sigma_1^2}{2}\!\sum_{k=1}^{d_1} \!\bbE\!\left[\exp\left(-\alpha\CE(\OX_t,\bbE\OY_t)\right)\! \left(\OX_t \!-\! \conspointx{\rho_{X,t}}\right)_k^2\!\cdot\!\left(\alpha^2\!\left(\partial_{x_k}\CE(\OX_t,\bbE\OY_t)\right)^2\!-\!\alpha\partial^2_{x_kx_k}\CE(\OX_t,\bbE\OY_t)\right)\right] dt\\
		&\;\; -\alpha \bbE\left[\exp\left(-\alpha\CE(\OX_t,\bbE\OY_t)\right)\nabla_y\CE(\OX_t,\bbE\OY_t) \cdot d\bbE\OY_t\right]
		=: \left(T_1+T_2+T_3\right)dt,
	\end{split}
	\end{equation*}
	where for the definition in the last step we exploited that the appearing stochastic integrals have expectation $0$ as a consequence of the regularity established in Theorem~\ref{theorem:well-posedness_mean-field} and Assumptions~\ref{ass:sec:convergence_1} and \ref{ass:sec:convergence_2}.
	Noticing that $\bbE\left[\exp\left(-\alpha\CE(\OX_t,\bbE\OY_t)\right) \left(\OX_t - \conspointx{\rho_{X,t}}\right)\right] = 0$ and $\nabla_x\CE(\conspointx{\rho_{X,t}},\bbE\OY_t)$ is deterministic, we obtain for $T_1$ the lower bound
	\begin{equation*} 
	\begin{split}
		T_1 
		&= \alpha\lambda_1 \bbE\left[\exp\left(-\alpha\CE(\OX_t,\bbE\OY_t)\right)\nabla_x\CE(\OX_t,\bbE\OY_t) \cdot \left(\OX_t - \conspointx{\rho_{X,t}}\right)\right] \\
		&= \alpha\lambda_1 \bbE\left[\exp\left(-\alpha\CE(\OX_t,\bbE\OY_t)\right)\left(\nabla_x\CE(\OX_t,\bbE\OY_t)-\nabla_x\CE(\conspointx{\rho_{X,t}},\bbE\OY_t)\right) \cdot \left(\OX_t - \conspointx{\rho_{X,t}}\right)\right] \\
		&\geq -\alpha\lambda_1 e^{-\alpha\underbar{\CE}(\bbE\overbarscript{Y}_t)}C_{\nabla^2\CE}\bbE\N{\OX_t - \conspointx{\rho_{X,t}}}_2^2,
	\end{split}
	\end{equation*}
	where we made use of the assumptions again.
	For $T_2$ it holds 
	\begin{equation*} 
	\begin{split}
		T_2
		&\geq -\alpha\frac{\sigma_1^2}{2}\sum_{k=1}^{d_1} \bbE\left[\exp\left(-\alpha\CE(\OX_t,\bbE\OY_t)\right) \left(\OX_t - \conspointx{\rho_{X,t}}\right)_k^2\partial^2_{x_kx_k}\CE(\OX_t,\bbE\OY_t)\right] dt\\
		&\geq - \alpha\frac{\sigma_1^2}{2}e^{-\alpha\underbar{\CE}(\bbE\overbarscript{Y}_t)}C_{\nabla^2\CE} \bbE\N{\OX_t - \conspointx{\rho_{X,t}}}_2^2.
	\end{split}
	\end{equation*}
	And, eventually, for $T_3$ we have the following bound from below
	\begin{equation*} 
	\begin{split}
		T_3
		&= \alpha\lambda_2 \bbE\left[\exp\left(-\alpha\CE(\OX_t,\bbE\OY_t)\right)\nabla_y\CE(\OX_t,\bbE\OY_t) \cdot \left(\bbE\OY_t - \conspointy{\rho_{Y,t}}\right)\right] \\
		&\geq -\alpha\lambda_2 e^{-\alpha\underbar{\CE}(\bbE\overbarscript{Y}_t)} C_{\nabla\CE} \N{\bbE\OY_t - \conspointy{\rho_{Y,t}}}_2,
	\end{split}
	\end{equation*}
	where we used the bounds on the gradient of $\CE$ required through Assumption~\ref{ass:sec:convergence_1} in the last step.
	Collecting the estimates for $T_1$, $T_2$ and $T_3$, and inserting them into the first equation gives
	\begin{equation*} 
	\begin{split}
		\frac{d}{dt}\widetilde\CM^X(t)
		&\!\geq\! -\alpha e^{-\alpha\underbar{\CE}(\bbE\overbarscript{Y}_t)} C_{\nabla^2\CE} \!\left(\lambda_1 \!+\! \frac{\sigma_1^2}{2}\right) \!\bbE\N{\OX_t \!-\! \conspointx{\rho_{X,t}}}_2^2
		\!-\!\alpha\lambda_2 e^{-\alpha\underbar{\CE}(\bbE\overbarscript{Y}_t)} C_{\nabla\CE} \N{\bbE\OY_t \!-\! \conspointy{\rho_{Y,t}}}_2\!.
	\end{split}
	\end{equation*}
	The two appearing norms can be bounded by recalling~\eqref{eq:lem:evolution_Var_6} and noticing that~\eqref{eq:lem:evolution_Var_5} gives 
	\begin{equation} \label{eq:lem:evolution_CM_7}
	\begin{split}
		\N{\bbE\OY_t-\conspointy{\rho_{Y,t}}}_2^2
		&\leq \frac{1}{\bbE\omegab(\bbE\OX_t,\OY_t)} e^{\beta\overbar{\CE}(\bbE\overbarscript{X}_t)} \int \! \N{\bbE\OY_t-y}_2^2 \,d\rho_{Y,t}(y)
		\leq \frac{\VarY(t)}{\CM^Y(t)}.
	\end{split}
	\end{equation}
	Inserting these two latter estimates allows to continue the former as desired as
	\begin{equation} \label{eq:lem:evolution_CM_8}
	\begin{split}
		\frac{d}{dt}\widetilde\CM^X(t)
		&\geq -4\alpha e^{-\alpha\underbar{\CE}(\bbE\overbarscript{Y}_t)} C_{\nabla^2\CE} \left(\lambda_1 + \frac{\sigma_1^2}{2}\right) \frac{\VarX(t)}{\CM^X(t)}-\alpha\lambda_2 e^{-\alpha\underbar{\CE}(\bbE\overbarscript{Y}_t)} C_{\nabla\CE} \frac{\sqrt{\VarY(t)}}{\sqrt{\CM^Y(t)}}.
	\end{split}
	\end{equation}
	The estimate for $\frac{d}{dt}\widetilde\CM^Y(t)$ can be obtained analogously.
\end{proof}

As mentioned already before we derive in the next lemma the time-evolutions of the functionals $\CM^X$ and $\CM^Y$ \revised{as defined in \eqref{eq:CM}.}
This is an immediate consequence of product rule and Lemma~\ref{lem:evolution_CM_aux}.

\begin{lemma} \label{lem:evolution_CM}
	Let $\VarX$ and $\VarY$ be as defined in~\eqref{eq:variances}, and $\CM^X$ and $\CM^Y$ as in~\eqref{eq:CM}.
	Then, under Assumptions~\ref{ass:sec:convergence_1} and \ref{ass:sec:convergence_2}, it \revised{holds}
	\begin{equation}
	\begin{split}
		\frac{d}{dt}\CM^X(t)
		&\geq -4\alpha C_{\nabla^2\CE} \left(\lambda_1 + \frac{\sigma_1^2}{2}\right) \frac{\VarX(t)}{\CM^X(t)} 
		-2\alpha\lambda_2 C_{\nabla\CE} \frac{\sqrt{\VarY(t)}}{\sqrt{\CM^Y(t)}}
	\end{split}
	\end{equation}
	as well as
	\begin{equation}
	\begin{split}
		\frac{d}{dt}\CM^Y(t)
		&\geq -4\beta C_{\nabla^2\CE} \left(\lambda_2 + \frac{\sigma_2^2}{2}\right) \frac{\VarY(t)}{\CM^Y(t)}
		-2\beta\lambda_1 C_{\nabla\CE} \frac{\sqrt{\VarX(t)}}{\sqrt{\CM^X(t)}}.
	\end{split}
	\end{equation}
\end{lemma}

\begin{proof}
	By product rule we have $\frac{d}{dt}\CM^X(t) = e^{\alpha\underbar{\CE}(\bbE\overbarscript{Y}_t)} \frac{d}{dt}\widetilde\CM^X(t) + \widetilde\CM^X(t)\,\frac{d}{dt}e^{\alpha\underbar{\CE}(\bbE\overbarscript{Y}_t)}$.
	While the first summand is controlled by recalling Lemma~\ref{lem:evolution_CM_aux}, for the second term we straightforwardly compute
	\begin{equation*}
	\begin{split}
		\frac{d}{dt}e^{\alpha\underbar{\CE}(\bbE\overbarscript{Y}_t)}
		&= \alpha e^{\alpha\underbar{\CE}(\bbE\overbarscript{Y}_t)} \nabla\underbar{\CE}(\bbE\overbarscript{Y}_t) \cdot \frac{d}{dt} \bbE\overbarscript{Y}_t
		\geq -\alpha\lambda_2 e^{\alpha\underbar{\CE}(\bbE\overbarscript{Y}_t)} C_{\nabla\CE} \N{\bbE\OY_t - \conspointy{\rho_{Y,t}}}_2,
	\end{split}
	\end{equation*}
	where we used the bounds on the gradient of $\underbar\CE$ required through Assumption~\ref{ass:sec:convergence_1} together with the regularity from Theorem~\ref{theorem:well-posedness_mean-field}.
	Recalling \eqref{eq:lem:evolution_CM_7} and putting everything together yields
	\begin{equation*}
	\begin{split}
		\frac{d}{dt}\CM^X(t)
		&\geq -4\alpha \revised{C_{\nabla^2\CE}} \left(\lambda_1 + \frac{\sigma_1^2}{2}\right) \frac{\VarX(t)}{\CM^X(t)}-\revised{\left(1+\CM^X(t)\right)}\alpha\lambda_2 C_{\nabla\CE} \frac{\sqrt{\VarY(t)}}{\sqrt{\CM^Y(t)}},
	\end{split}
	\end{equation*}
	\revised{which gives the claim after noting that $\CM^X(t)\leq1$.}
	The proceeding for $\frac{d}{dt}\CM^Y(t)$ is identical.
\end{proof}

\subsection{Time-Evolution of the Functionals $\CM^X_*$ and $\CM^Y_*$ \revised{from \eqref{eq:CM_star}}} \label{subsec:evolution_CM*}

Similarly to the preceding sections we study the time-evolution of two functionals $\CM^X_*$ and $\CM^Y_*$ as defined in \eqref{eq:CM_star}, which aids to prove properties of the limit point of the mean-field dynamics~\eqref{eq:saddle_point_dynamics_macro}.

\begin{lemma} \label{lem:evolution_CM*}
	Let $\VarX$ and $\VarY$ be as defined in~\eqref{eq:variances}, $\CM^X$ and $\CM^Y$ as in~\eqref{eq:CM}, and $\CM^X_*$ and $\CM^Y_*$ as in~\eqref{eq:CM_star}.
	Then, under Assumptions~\ref{ass:sec:convergence_1} and \ref{ass:sec:convergence_2}, it \revised{holds}
	\begin{equation}
	\begin{split}
		\frac{d}{dt}\CM^X_*(t)
		&\geq -4\alpha\lambda_1 e^{-\alpha\underbar{\CE}(y^*)}C_{\nabla\CE}\frac{\sqrt{\VarX(t)}}{\sqrt{\CM^X(t)}}- 2\alpha\sigma_1^2e^{-\alpha\underbar{\CE}(y^*)}C_{\nabla^2\CE} \frac{\VarX(t)}{\CM^X(t)}
	\end{split}
	\end{equation}
	as well as
	\begin{equation}
	\begin{split}
		\frac{d}{dt}\CM^Y_*(t)
		&\geq -4\beta\lambda_2 e^{\beta\overbar{\CE}(x^*)}C_{\nabla\CE}\frac{\sqrt{\VarY(t)}}{\sqrt{\CM^Y(t)}}- 2\beta\sigma_2^2e^{\beta\overbar{\CE}(x^*)}C_{\nabla^2\CE} \frac{\VarY(t)}{\CM^Y(t)}.
	\end{split}
	\end{equation}
\end{lemma}

\begin{proof}
	With It\^o's formula and chain rule we first note that 
	\begin{equation} \label{eq:lem:evolution_CM*_1}
	\begin{split}
		d\CM^X_*(t) &= 
		-\alpha \bbE\left[\exp\left(-\alpha\CE(\OX_t,y^*)\right)\nabla_x\CE(\OX_t,y^*) \cdot d\OX_t\right]+\frac{\sigma_1^2}{2}\sum_{k=1}^{d_1} \bbE\left[\exp\left(-\alpha\CE(\OX_t,y^*)\right)\cdot\right.\\
		&\quad\quad\left.\cdot\left(\OX_t \!-\! \conspointx{\rho_{X,t}}\right)_k^2\left(\alpha^2\left(\partial_{x_k}\CE(\OX_t,y^*)\right)^2\!-\!\alpha\partial^2_{x_kx_k}\CE(\OX_t,y^*)\right)\right] dt
		=: \left(T_1+T_2\right)dt,
	\end{split}
	\end{equation}
	where for the definition in the last step we again exploited that the appearing stochastic integral has expectation $0$ as a consequence of the assumptions.
	For $T_1$ we have the lower bound
	\begin{equation*} 
	\begin{split}
		T_1 
		&\!\geq\! -\alpha\lambda_1 e^{-\alpha\underbar{\CE}(y^*)} \bbE\left[\N{\nabla_x\CE(\OX_t,y^*)}_2 \N{\OX_t \!-\! \conspointx{\rho_{X,t}}}_2\right] 
		\!\geq\! -\alpha\lambda_1 e^{-\alpha\underbar{\CE}(y^*)}C_{\nabla\CE}\sqrt{\bbE\N{\OX_t \!-\! \conspointx{\rho_{X,t}}}_2^2}.
	\end{split}
	\end{equation*}
	For $T_2$ it holds
	\begin{equation*} 
	\begin{split}
		T_2
		&\geq -\alpha\frac{\sigma_1^2}{2}\sum_{k=1}^{d_1} \bbE\left[\exp\left(-\alpha\CE(\OX_t,y^*)\right) \left(\OX_t - \conspointx{\rho_{X,t}}\right)_k^2\partial^2_{x_kx_k}\CE(\OX_t,y^*)\right] dt\\
		&\geq - \alpha\frac{\sigma_1^2}{2}e^{-\alpha\underbar{\CE}(y^*)}C_{\nabla^2\CE}\bbE\N{\OX_t - \conspointx{\rho_{X,t}}}_2^2.
	\end{split}
	\end{equation*} 
	Collecting the two former estimates for the terms $T_1$ and $T_2$, and inserting them into~\eqref{eq:lem:evolution_CM*_1} gives
	\begin{equation} \label{eq:lem:evolution_CM*_5}
	\begin{split}
		\frac{d}{dt}\CM^X_*(t)
		&\!\geq\! -\alpha\lambda_1 e^{-\alpha\underbar{\CE}(y^*)}C_{\nabla\CE}\sqrt{\bbE\N{\OX_t \!-\! \conspointx{\rho_{X,t}}}_2^2} \!-\!  \alpha\frac{\sigma_1^2}{2}e^{-\alpha\underbar{\CE}(y^*)}C_{\nabla^2\CE}\bbE\N{\OX_t \!-\! \conspointx{\rho_{X,t}}}_2^2,
	\end{split}
	\end{equation}
	where the last expression can be bounded by employing~\eqref{eq:lem:evolution_Var_6}.
	The estimate for $\frac{d}{dt}\CM^Y_*(t)$ can be obtained analogously.
\end{proof}

\subsection{Proof of Theorem~\ref{theorem:convergece}} \label{subsec:proof_main_theorem}


\begin{proof}[Proof of Theorem~\ref{theorem:convergece}]
	\noindent\textbf{Step 1a:}
	Let us define the time horizon
	\begin{equation} \label{eq:thm:convergece_1}
		T := \inf\left\{t\geq0 : \CM^X(t)<\frac{1}{2}\CM^X(0) \text{ or }\CM^Y(t)<\frac{1}{2}\CM^Y(0) \right\}
		\quad
		\text{with } \inf\emptyset=\infty,
	\end{equation}
	\revised{where $\CM^X$ and $\CM^Y$ are as defined in \eqref{eq:CM}.}
	Obviously, by continuity, $T>0$.
	We claim that $T=\infty$, which is shown by contradiction in what follows.
	Therefore, let us assume $T<\infty$.
	Then, as a consequence of the definition of the time horizon $T$, the prefactors of $\VarX(t)$ and $\VarY(t)$ in Lemma~\ref{lem:evolution_Var} are upper bounded by $-\mu_1$ and $-\mu_2$, respectively, for all $t\in[0,T]$.
	Consequently, Lemma~\ref{lem:evolution_Var} permits the upper bounds
	\begin{equation} \label{eq:thm:convergece_2}
		\frac{d}{dt}\VarX(t) \leq -\mu_1 \VarX(t)
		\quad\text{and}\quad
		\frac{d}{dt}\VarY(t) \leq -\mu_2 \VarY(t)
	\end{equation}
	for the time-evolution of the functionals~$\VarX$ and~$\VarY$.
	The negativity of the rate is ensured by the well-preparedness condition~\ref{asm:well_preparedness_P1}. 
	An application of Gr\"onwall's inequality gives
	\begin{equation} \label{eq:thm:convergece_3}
		\VarX(t) \leq \VarX(0)e^{-\mu_1t}
		\quad\text{and}\quad
		\VarY(t) \leq \VarY(0)e^{-\mu_2t}.
	\end{equation}
	Let us now derive the contradiction.
	It follows from Lemma~\ref{lem:evolution_CM} \revised{for $\CM^X$ and $\CM^Y$ from \eqref{eq:CM}} that
	\begin{equation} \label{eq:thm:convergece_4}
	\begin{split}
		\frac{d}{dt}\CM^X(t)
		&\geq -8\alpha C_{\nabla^2\CE} \left(\lambda_1 + \frac{\sigma_1^2}{2}\right) \frac{\VarX(0)e^{-\mu_1t}}{\CM^X(0)} 
		-2\sqrt{2}\alpha\lambda_2C_{\nabla\CE} \frac{\sqrt{\VarY(0)}e^{-\mu_2t/2}}{\sqrt{\CM^Y(0)}},\\
		\frac{d}{dt}\CM^Y(t)
		&\geq -8\beta C_{\nabla^2\CE} \left(\lambda_2 + \frac{\sigma_2^2}{2}\right) \frac{\VarY(0)e^{-\mu_2t}}{\CM^Y(0)}
		-2\sqrt{2}\beta\lambda_1C_{\nabla\CE} \frac{\sqrt{\VarX(0)}e^{-\mu_1t/2}}{\sqrt{\CM^X(0)}}
	\end{split}
	\end{equation}
	where we used the formerly derived~\eqref{eq:thm:convergece_3} as well as that $\CM^X(t)\geq\CM^X(0)/2$ and $\CM^Y(t)\geq\CM^Y(0)/2$ for all $t\in[0,T]$ by definition of $T$.
	Integrating~\eqref{eq:thm:convergece_4} and employing the well-preparedness condition~\ref{asm:well_preparedness_P2} shows for all $t\in[0,T]$
	\begin{equation*} 
	\begin{split}
		\CM^X(t)
		&\geq \CM^X(0)
		-8\alpha C_{\nabla^2\CE} \left(\lambda_1 + \frac{\sigma_1^2}{2}\right) \frac{\VarX(0)}{\mu_1\CM^X(0)}
		-4\sqrt{2}\alpha\lambda_2 C_{\nabla\CE} \frac{\sqrt{\VarY(0)}}{\mu_2\sqrt{\CM^Y(0)}}
		\geq\frac{3}{4}\CM^X(0),\\
		\CM^Y(t)
		&\geq \CM^Y(0)
		-8\beta C_{\nabla^2\CE} \left(\lambda_2 + \frac{\sigma_2^2}{2}\right) \frac{\VarY(0)}{\mu_2\CM^Y(0)}
		-4\sqrt{2}\beta\lambda_1 C_{\nabla\CE} \frac{\sqrt{\VarX(0)}}{\mu_1\sqrt{\CM^X(0)}}
		\geq\frac{3}{4}\CM^Y(0).
	\end{split}
	\end{equation*}
	This entails that there exists $\delta>0$ such that $\CM^X(t)\geq\CM^X(0)/2$ and $\CM^Y(t)\geq\CM^Y(0)/2$ hold for all $t\in[T,T+\delta]$ as well, contradicting the definition of $T$ and therefore showing $T=\infty$.
	Consequently \eqref{eq:thm:convergece_3} as well as 
	\begin{equation} \label{eq:thm:convergece_8}
		\CM^X(t)\geq\frac{1}{2}\CM^X(0)
		\quad\text{and}\quad
		\CM^Y(t)\geq\frac{1}{2}\CM^Y(0)
	\end{equation}
	hold for all $t\geq0$, which proves~\eqref{eq:theorem:convergece:3}.
	
	\noindent\textbf{Step 1b:}
	With Jensen's inequality and by making use of the bounds~\eqref{eq:lem:evolution_Var_6} and~\eqref{eq:lem:evolution_Var_7} combined with \eqref{eq:thm:convergece_3} and~\eqref{eq:thm:convergece_8} we further observe that
	\begin{equation*} 
	\begin{split}
		\N{\frac{d}{dt}\bbE \OX_t}_2
		&\leq \lambda_1\bbE\N{\OX_t-\conspointx{\rho_{X,t}}}_2
		\leq 2\lambda_1 \frac{\sqrt{\VarX(0)}e^{-\mu_1t/2}}{\sqrt{\CM^X(0)}}
		\rightarrow 0 \quad\text{ as } t\rightarrow\infty,\\
		\N{\frac{d}{dt}\bbE \OY_t}_2
		&\leq \lambda_2\bbE\N{\OY_t-\conspointy{\rho_{Y,t}}}_2
		\leq 2\lambda_2 \frac{\sqrt{\VarY(0)}e^{-\mu_2t/2}}{\sqrt{\CM^Y(0)}}
		\rightarrow 0 \quad\text{ as } t\rightarrow\infty.
	\end{split}
	\end{equation*}
	We therefore have $\left(\bbE \OX_t, \bbE \OY_t\right)\rightarrow\left(\widetilde{x},\widetilde{y}\right)$ for \textit{some} $\left(\widetilde{x},\widetilde{y}\right)\in\bbR^{d_1+d_2}$.
	In fact, following from~\eqref{eq:thm:convergece_3}, $\left( \OX_t, \OY_t\right)\rightarrow\left(\widetilde{x},\widetilde{y}\right)$  and $\big(\conspointx{\rho_{X,t}},\conspointy{\rho_{Y,t}}\big)\rightarrow\left(\widetilde{x},\widetilde{y}\right)$ in $L^2$ thanks to \eqref{eq:lem:evolution_Var_6} and~\eqref{eq:lem:evolution_Var_7}.
	This shows~\eqref{eq:theorem:convergece:4}.
	
	\noindent\textbf{Step 2a:}
	It remains to \revised{verify \eqref{eq:theorem:convergece:5} for the point~$\left(\widetilde{x},\widetilde{y}\right)$.}
	With similar arguments as in Step~1a let us first derive analogous statements as in~\eqref{eq:thm:convergece_8} for $\widetilde\CM^X$ and $\widetilde\CM^Y$ \revised{as defined in \eqref{eq:CM_tilde}} as well as $\CM^X_*$ and $\CM^Y_*$ \revised{as defined in \eqref{eq:CM_star},} respectively.
	To do so, we first notice that $(\bbE \OX_t,\bbE \OY_t)$ is continuous and since it converges to $(\widetilde{x},\widetilde{y})$ as $t\rightarrow\infty$, there exists $M>0$, potentially depending on $(\widetilde{x},\widetilde{y})$, such that $\N{\bbE \OX_t}_2+\N{\bbE \OY_t}_2 \leq M$ for all $t\geq0$.
	Since moreover $\underbar{\CE}$ and $\overbar{\CE}$ are continuous, there exists $\CE_M>0$ such that $-\CE_M\leq\underbar{\CE}(\bbE \OY_t) \leq \overbar{\CE}(\bbE \OX_t) \leq \CE_M$ for all $t>0$.
	Utilizing this together with~\eqref{eq:thm:convergece_3} and~\eqref{eq:thm:convergece_8} we derive from Lemma~\ref{lem:evolution_CM_aux} \revised{for $\widetilde\CM^X$ and $\widetilde\CM^Y$ from \eqref{eq:CM_tilde}} that
	\begin{equation} \label{eq:thm:convergece_9a}
	\begin{split}
		\frac{d}{dt}\widetilde\CM^X(t)
		&\geq -8\alpha e^{\alpha\CE_M} C_{\nabla^2\CE} \left(\lambda_1 \!+\! \frac{\sigma_1^2}{2}\right) \frac{\VarX(0)e^{-\mu_1t}}{\CM^X(0)} 
		\!-\!\sqrt{2}\alpha\lambda_2 e^{\alpha\CE_M} C_{\nabla\CE} \frac{\sqrt{\VarY(0)}e^{-\mu_2t/2}}{\sqrt{\CM^Y(0)}},\\
		\frac{d}{dt}\widetilde\CM^Y(t)
		&\geq -8\beta e^{\beta\CE_M} C_{\nabla^2\CE} \left(\lambda_2 \!+\! \frac{\sigma_2^2}{2}\right) \frac{\VarY(0)e^{-\mu_2t}}{\CM^Y(0)}
		\!-\!\sqrt{2}\beta\lambda_1 e^{\beta\CE_M} C_{\nabla\CE} \frac{\sqrt{\VarX(0)}e^{-\mu_1t/2}}{\sqrt{\CM^X(0)}}.
	\end{split}
	\end{equation}
	Analogously, by using~\eqref{eq:thm:convergece_3} and~\eqref{eq:thm:convergece_8} it follows directly from Lemma~\ref{lem:evolution_CM*} \revised{for $\CM^X_*$ and $\CM^Y_*$ from \eqref{eq:CM_star}} that
	\begin{equation} \label{eq:thm:convergece_9b}
	\begin{split}
		\frac{d}{dt}\CM^X_*(t)
		&\geq -8\alpha\lambda_1 e^{-\alpha\underbar{\CE}(y^*)}C_{\nabla\CE}\frac{\sqrt{\VarX(0)}e^{-\mu_1t/2}}{\sqrt{\CM^X(0)}}
		- 4\alpha\sigma_1^2e^{-\alpha\underbar{\CE}(y^*)}C_{\nabla^2\CE} \frac{\VarX(0)e^{-\mu_1t}}{\CM^X(0)},\\
		\frac{d}{dt}\CM^Y_*(t)
		&\geq -8\beta\lambda_2 e^{\beta\overbar{\CE}(x^*)}C_{\nabla\CE}\frac{\sqrt{\VarY(0)}e^{-\mu_2t/2}}{\sqrt{\CM^Y(0)}}
		- 4\beta\sigma_2^2e^{\beta\overbar{\CE}(x^*)}C_{\nabla^2\CE} \frac{\VarY(0)e^{-\mu_2t}}{\CM^Y(0)}.
	\end{split}
	\end{equation}
	Integrating \eqref{eq:thm:convergece_9a} and employing the well-preparedness condition~\ref{asm:well_preparedness_P2} shows for all $t\geq0$ that
	\begin{equation*} 
	\begin{split}
		\widetilde\CM^X(t)
		&\!\geq\! \widetilde\CM^X(0)
		\!-\!8\alpha e^{\alpha\CE_M} C_{\nabla^2\CE} \!\left(\lambda_1 \!+\! \frac{\sigma_1^2}{2}\right) \frac{\VarX(0)}{\mu_1\CM^X(0)} \!-\!2\sqrt{2}\alpha\lambda_2 e^{\alpha\CE_M} C_{\nabla\CE} \frac{\sqrt{\VarY(0)}}{\mu_2\sqrt{\CM^Y(0)}} \!\geq\! \frac{3}{4}\widetilde\CM^X(0),\\
		\widetilde\CM^Y(t)
		&\!\geq\! \widetilde\CM^Y(0)
		\!-\!8\beta e^{\beta\CE_M} C_{\nabla^2\CE} \!\left(\lambda_2 \!+\! \frac{\sigma_2^2}{2}\right) \frac{\VarY(0)}{\mu_2\CM^Y(0)} \!-\!2\sqrt{2}\beta\lambda_1 e^{\beta\CE_M} C_{\nabla\CE} \frac{\sqrt{\VarX(0)}}{\mu_1\sqrt{\CM^X(0)}} \!\geq\! \frac{3}{4}\widetilde\CM^Y(0).
	\end{split}
	\end{equation*}
	Analogously, using \eqref{eq:thm:convergece_9b} together with~\ref{asm:well_preparedness_P2} shows for all $t\geq0$ that
	\begin{equation*} 
	\begin{split}
		\CM^X_*(t)
		&\!\geq\! \CM^X_*(0)
		\!-\!16\alpha\lambda_1 e^{-\alpha\underbar{\CE}(y^*)}C_{\nabla\CE}\frac{\sqrt{\VarX(0)}}{\mu_1\sqrt{\CM^X(0)}} \!-\! 4\alpha\sigma_1^2e^{-\alpha\underbar{\CE}(y^*)}C_{\nabla^2\CE} \frac{\VarX(0)}{\mu_1\CM^X(0)} \!\geq\! \frac{3}{4}\CM^X_*(0),\\
		\CM^Y_*(t)
		&\!\geq\! \CM^Y_*(0)
		\!-\!16\beta\lambda_2 e^{\beta\overbar{\CE}(x^*)}C_{\nabla\CE}\frac{\sqrt{\VarY(0)}}{\mu_2\sqrt{\CM^Y(0)}} \!-\! 4\beta\sigma_2^2e^{\beta\overbar{\CE}(x^*)}C_{\nabla^2\CE} \frac{\VarY(0)}{\mu_2\CM^Y(0)} \!\geq\! \frac{3}{4}\CM^Y_*(0).
	\end{split}
	\end{equation*}
	Thus, in particular it holds for all $t\geq0$ 
	\begin{equation} \label{eq:thm:convergece_13a}
		\widetilde\CM^X(t)\geq\frac{1}{2}\widetilde\CM^X(0)
		\quad\text{and}\quad
		\widetilde\CM^Y(t)\geq\frac{1}{2}\widetilde\CM^Y(0)
	\end{equation}
	as well as 
	\begin{equation} \label{eq:thm:convergece_13b}
		\CM^X_*(t)\geq\frac{1}{2}\CM^X_*(0)
		\quad\text{and}\quad
		\CM^Y_*(t)\geq\frac{1}{2}\CM^Y_*(0).
	\end{equation}
	
	\noindent\textbf{Step 2b:}
\revised{By Chebyshev's inequality, for each $\delta>0$ it holds that 
	\begin{align*}
		\rho_t(\{\N{(x-\widetilde{x},y-\widetilde{y})}_2\geq \delta\})
		&\leq \frac{2}{\delta^2}\big(\VarX(t)+\VarY(t)+\N{\bbE \OX_t-\widetilde{x}}_2^2+\N{\bbE \OY_t-\widetilde{y}}_2^2\big)
		\\
		&\rightarrow 0, \text{ as } t\rightarrow \infty.
		\end{align*}
	Thus,  the pair $(\OX_t,\OY_t)$ converges to $(\widetilde{x},\widetilde{y})$ in probability as $t$ tends to infinity.
Recall the convergence $\left(\bbE \OX_t, \bbE \OY_t\right)\rightarrow\left(\widetilde{x},\widetilde{y}\right)$, the continuity of $\CE$, and the fact that for all $t\geq 0$,
$$\exp(-\alpha \CE(\OX_t,\bbE \OY_t))\leq \exp(\alpha \CE_M), \quad \text{a.s.}$$
By the dominated convergence theorem, one can pass to the limit in $t$ to obtain
  $\lim_{t\rightarrow \infty}\widetilde\CM^X(t) = \exp\left(-\alpha\CE(\widetilde{x},\widetilde{y})\right)$. Analogously, one may get $\widetilde\CM^Y(t) \rightarrow \exp\left(\beta\CE(\widetilde{x},\widetilde{y})\right) $
	as $t\rightarrow\infty$.}
	Using this when taking the limit $t\rightarrow\infty$ in the bounds~\eqref{eq:thm:convergece_13a} after applying the logarithm and multiplying both sides with $-1/\alpha$ and $1/\beta$, respectively, we obtain
	\begin{equation} \label{eq:thm:convergece_17}
	\begin{split}
		\CE(\widetilde{x},\widetilde{y}) 
		&= \lim_{t\rightarrow\infty} \left(-\frac{1}{\alpha}\log\widetilde\CM^X(t)\right)
		\leq \frac{1}{\alpha}\log2 - \frac{1}{\alpha}\log\widetilde\CM^X(0),\\
		\CE(\widetilde{x},\widetilde{y})
		&= \lim_{t\rightarrow\infty} \left(\frac{1}{\beta}\log\widetilde\CM^Y(t)\right) 
		\revised{\geq -\frac{1}{\beta}\log2 + \frac{1}{\beta}\log\widetilde\CM^Y(0)}.
	\end{split}
	\end{equation}
	\revised{Due to the first set of well-preparedness conditions from \ref{asm:well_preparedness_P1'}, the Laplace principle in form of Lemmas~\ref{lmx} and \ref{lmy}} when choosing $\mu^\alpha$ as the law of the initial data $\OX_0$ and $\mu^\beta$ as the law of $\OY_0$, now allows to choose $\alpha\geq(2\log 2)/\varepsilon$ and $\beta\geq(2\log 2)/\varepsilon$ large enough such that for given $\varepsilon>0$ it moreover \revised{holds}
	\begin{equation} \label{eq:thm:convergece_19}
	\begin{split}
		-\frac{1}{\alpha}\log\widetilde\CM^X(0)-\min_{x\in\bbR^{d_1}}\CE(x,\bbE\OY_0)
		&= -\frac{1}{\alpha}\log\bbE\exp\left(-\alpha\CE(\OX_0,\bbE\OY_0)\right)-\min_{x\in\bbR^{d_1}}\CE(x,\bbE\OY_0)
		\leq \varepsilon/2, \\
		-\frac{1}{\beta}\log\widetilde\CM^Y(0)+\max_{y\in\bbR^{d_2}}\CE(\bbE\OX_0,y)
		&=-\frac{1}{\beta}\log\bbE\exp\left(\beta\CE(\bbE\OX_0,\OY_0)\right)+\max_{y\in\bbR^{d_2}}\CE(\bbE\OX_0,y)
		\leq \varepsilon/2.
	\end{split}	
	\end{equation}
	Notice here that we well-prepare $\alpha$ and $\beta$ simultaneously with the initial data $(\OX_0,\OY_0)$ (therewith $(\OX_0,\OY_0)$ depends on $\alpha$, $\beta$).
	However due to the well-preparedness conditions~\ref{asm:well_preparedness_P1'},  $\alpha$ and $\beta$ can still be taken sufficiently large as ensured in Lemmas~\ref{lmx} and \ref{lmy}.

	Such choices of parameters in Equation~\eqref{eq:thm:convergece_17} immediately give
	\begin{equation} \label{eq:thm:convergece_21}
		\CE(\widetilde{x},\widetilde{y}) \leq \min_{x\in\bbR^{d_1}}\CE(x,\bbE\OY_0) + \varepsilon
		\quad\text{and}\quad
		\revised{\CE(\widetilde{x},\widetilde{y}) \geq \max_{y\in\bbR^{d_2}}\CE(\bbE\OX_0,y) - \varepsilon}
	\end{equation}
	and consequently 
	\begin{equation} \label{eq:thm:convergece_22}
		\CE(\widetilde{x},\widetilde{y}) \leq \min_{x\in\bbR^{d_1}}\max_{y\in\bbR^{d_2}}\CE(x,y) + \varepsilon
		\quad\text{and}\quad
		\revised{\CE(\widetilde{x},\widetilde{y}) \geq \max_{y\in\bbR^{d_2}}\min_{x\in\bbR^{d_1}}\CE(x,y) -\varepsilon},
	\end{equation}
	which proves the first part of \eqref{eq:theorem:convergece:5}.
	Secondly, following an analogous argumentation for $\CM^X_*$ and $\CM^Y_*$ \revised{as defined in \eqref{eq:CM_star},} we obtain the remainder of~\eqref{eq:theorem:convergece:5}.
	More precisely, we first note that, as $t\rightarrow\infty$,
	\begin{equation} \label{eq:thm:convergece_23}
		\CM^X_*(t) \rightarrow \exp\left(-\alpha\CE(\widetilde{x},y^*)\right)
		\quad\text{and}\quad
		\CM^Y_*(t) \rightarrow \exp\left(\beta\CE(x^*,\widetilde{y})\right). 
	\end{equation}
	Taking now the limit $t\rightarrow\infty$ in~\eqref{eq:thm:convergece_13b} after suitable algebraic manipulations, we obtain
	\begin{equation} \label{eq:thm:convergece_24}
	\begin{split}
		\CE(\widetilde{x},y^*) 
		&= \lim_{t\rightarrow\infty} \left(-\frac{1}{\alpha}\log\CM^X_*(t)\right)
		\leq \frac{1}{\alpha}\log2 - \frac{1}{\alpha}\log\CM^X_*(0),\\
		\CE(x^*,\widetilde{y})
		&= \lim_{t\rightarrow\infty} \left(\frac{1}{\beta}\log\CM^Y_*(t)\right) 
		\geq -\frac{1}{\beta}\log2 + \frac{1}{\beta}\log\CM^Y_*(0).
	\end{split}
	\end{equation}
	A potentially larger choice of $\alpha$ and $\beta$ allows (\revised{again by the Laplace principle in form of Lemmas~\ref{lmx} and \ref{lmy}, which applies due to the second set of well-preparedness conditions from \ref{asm:well_preparedness_P1'})} to guarantee 
	\begin{equation*} 
	\begin{split}
		-\frac{1}{\alpha}\log\CM^X_*(0)-\min_{x\in\bbR^{d_1}}\CE(x,y^*)
		= -\frac{1}{\alpha}\log\bbE\exp\left(-\alpha\CE(\OX_0,y^*)\right)-\min_{x\in\bbR^{d_1}}\CE(x,y^*)
		&\leq \varepsilon/2,\\
		-\frac{1}{\beta}\log\CM^Y_*(0)+\max_{y\in\bbR^{d_2}}\CE(x^*,y)
		= -\frac{1}{\beta}\log\bbE\exp\left(\beta\CE(x^*,\OY_0)\right)+\max_{y\in\bbR^{d_2}}\CE(x^*,y)
		&\leq \varepsilon/2
	\end{split}
	\end{equation*}
	for the specified $\varepsilon$.
	Such choices of parameters in Equation~\eqref{eq:thm:convergece_24} immediately give
	\begin{equation} \label{eq:thm:convergece_28}
		\CE(\widetilde{x},y^*) \leq \min_{x\in\bbR^{d_1}}\CE(x,y^*) + \varepsilon
		\quad\text{and}\quad
		\revised{\CE(x^*,\widetilde{y}) \geq \max_{y\in\bbR^{d_2}}\CE(x^*,y) - \varepsilon}\,,
	\end{equation}
which completes the proof of \eqref{eq:theorem:convergece:5}.
	
	\noindent\textbf{Step 3:} Finally, under the inverse continuity property~\ref{ass:sec:ICP} and making use of what we just proved, we additionally obtain~\eqref{eq:theorem:convergece:7}, which concludes the proof.
\end{proof}

\section{\revised{Implementation of CBO-SP and Numerical Experiments}} \label{sec:numerics}

\revised{\subsection{Numerical Algorithm and Implementation}
\label{sec:implementation}
In order to implement and run CBO-SP on a computer, we first fix a discrete time step size $\Delta t$ as well as a number of iterations~$K$ or define any other suitable stopping criterion.
Then, by discretizing the interacting particle system~\eqref{eq:saddle_point_dynamics_micro} via an Euler-Maruyama time discretization~\cite{higham2001algorithmic,platen1999introduction} as
\begin{subequations} \label{eq:saddle_point_dynamics_micro_discrete}
\begin{align}
	\widehat{X}_{k+1}^i
		&= \widehat{X}_{k}^i -\lambda_1\Delta t\left(\widehat{X}_k^i - \conspointx{\empmeasureX{k}}\right)
		+ \sigma_1 D\!\left(\widehat{X}_k^i - \conspointx{\empmeasureX{k}}\right) B_k^{X,i},
		\label{eq:saddle_point_dynamics_micro_discrete_X} \\
	\widehat{Y}_{k+1}^i
		&= \widehat{Y}_{k}^i -\lambda_2\Delta t\left(\widehat{Y}_k^i - \conspointy{\empmeasureY{k}}\right)
		+ \sigma_2 D\!\left(\widehat{Y}_k^i - \conspointy{\empmeasureY{k}}\right) B_k^{Y,i},
	\label{eq:saddle_point_dynamics_micro_discrete_Y}
\end{align}
\end{subequations}
where $\empmeasureX{k}$ and $\empmeasureY{k}$ denote the empirical averages of the iterates $(\widehat{X}_k^i)_{i=1,\dots,N_1}$ and $(\widehat{Y}_k^i)_{i=1,\dots,N_2}$ and where
\begin{subequations}
\begin{align} 
	\widehat{x}_\alpha^Y(\empmeasureX{k})
	&= \!\int \! x \,\frac{\omegaa\big(x,\int \!y\,d\empmeasureY{k}(y)\big)}{\Nbig{\omegaa\big(\,\cdot\,,\int \!y\,d\empmeasureY{k}(y)\big)}_{L_1(\empmeasureX{k})}}\,d\empmeasureX{k}(x), \label{eq:conspoint_x_discrete}\\
	\widehat{y}_\beta^X(\empmeasureY{k})
	&= \!\int \! y \,\frac{\omegab\big(\int \!x\,d\empmeasureX{k+1}(x),y\big)}{\Nbig{\omegab\big(\int \!x\,d\empmeasureX{k+1}(x),\,\cdot\,\big)}_{L_1(\empmeasureY{k})}}\,d\empmeasureY{k}(y),\label{eq:conspoint_y_discrete}
\end{align}
\end{subequations}
we obtain the implementable iterative scheme, which is used in the formulation of Algorithm~\ref{algorithm:CBOSP}.
Moreover, $\big((B_k^{X,i})_{k=1,\dots,K}\big)_{i=1,\dots,\revised{N_1}}$ and $\big((B_k^{Y,i})_{k=1,\dots,K}\big)_{i=1,\dots,\revised{N_2}}$ in \eqref{eq:saddle_point_dynamics_micro_discrete} are independent Gaussian vectors in $\bbR^{d_1}$ and $\bbR^{d_2}$, respectively, with covariance matrix $\Delta t \Id$.
Note that in Equation~\eqref{eq:conspoint_y_discrete} we could also use the old iterates~$\empmeasureX{k}$ instead of the new ones~$\empmeasureX{k+1}$ for the computation.}

\begin{algorithm}[!ht]
\caption{CBO-SP}
\revised{\begin{algorithmic}[1]
		\floatname{algorithm}{Procedure}
		\renewcommand{\algorithmicrequire}{\textbf{Input:}}
		\renewcommand{\algorithmicensure}{\textbf{Output:}}
		\renewcommand{\algorithmicloop}{\textbf{while }}
		\renewcommand\algorithmicdo{}
		\renewcommand\algorithmicthen{}
		
	\Require{Objective~$\CE$, discrete time step size~$\Delta t$, number of iterates~$K$, parameters~$\lambda_1,\lambda_2,\sigma_1,\sigma_2,\alpha,\beta$, number of particles~$N_1$ and $N_2$, initialization~$\rho_0$}
	\Ensure{Approximation~$\big(\widehat{x}_\alpha^Y(\empmeasureX{k}), \widehat{y}_\beta^X(\empmeasureY{k})\big)$ of the saddle point~$(\saddlepoint)$ of $\CE$}
	\State
		Generate the particles' initial positions $(X^i_0)_{i=1,\dots,N_1}$ and $(Y^i_0)_{i=1,\dots,N_2}$ according to the initial laws~$\rho_{X,0}$ and $\rho_{Y,0}$, respectively. Set $k=0$.
	\Loop{$k\leq K$ or stopping criterion not fulfilled}
		\State \parbox[t]{\dimexpr\linewidth-4.4em}{Compute the component $\widehat{x}_\alpha^Y(\empmeasureX{k})$ of the consensus point according to~\eqref{eq:conspoint_x_discrete}.\strut}
		\State \parbox[t]{\dimexpr\linewidth-4.4em}{Update the $X$-positions by computing $\big(\widehat{X}^i_{k+1}\big)_{i=1,\dots,N_1}$ according to~\eqref{eq:saddle_point_dynamics_micro_discrete_X}.\strut}
		\State \parbox[t]{\dimexpr\linewidth-4.4em}{Compute the component $ \widehat{y}_\beta^X(\empmeasureY{k})$ of the consensus point according to~\eqref{eq:conspoint_y_discrete}.\strut}
		\State \parbox[t]{\dimexpr\linewidth-4.4em}{Update the $Y$-positions by computing $\big(\widehat{Y}^i_{k+1}\big)_{i=1,\dots,N_2}$ according to~\eqref{eq:saddle_point_dynamics_micro_discrete_Y}.\strut}
		\State \parbox[t]{\dimexpr\linewidth-1.4em}{Check the stopping criterion and \textbf{break} if fulfilled. If not, continue and set $k=k+1$.\strut}
	\EndLoop
	\State Compute consensus point $\big(\widehat{x}_\alpha^Y(\empmeasureX{k}), \widehat{y}_\beta^X(\empmeasureY{k})\big)$ as final approximation to saddle point~$(\saddlepoint)$.
\end{algorithmic}}
	\label{algorithm:CBOSP}
\end{algorithm}

\subsection{Illustrative Numerical Experiments for CBO-SP}
\begin{figure}[htp!]
	\def\heightt{0.11}
	\def\widthh{0.24}
	\def\horzspace{0.02em}
	\centering
	\begin{subfigure}[b]{\widthh\textwidth}
        \centering
        \renewcommand{\thesubfigure}{a}
        \includegraphics[trim=85 250 79 250,clip,height=\heightt\textheight]{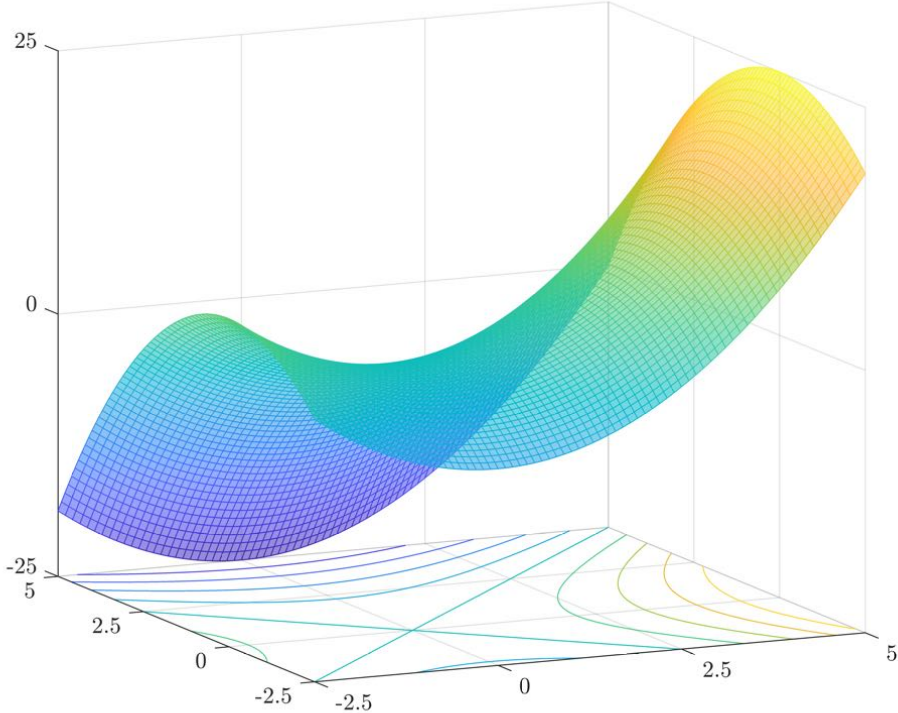}
        \caption{\footnotesize The objective function $\CE(x,y) = x^2-y^2$}
    \end{subfigure}~\hspace{\horzspace}~
	\begin{subfigure}[b]{\widthh\textwidth}
        \centering
        \renewcommand{\thesubfigure}{b}
        \includegraphics[trim=85 250 79 250,clip,height=\heightt\textheight]{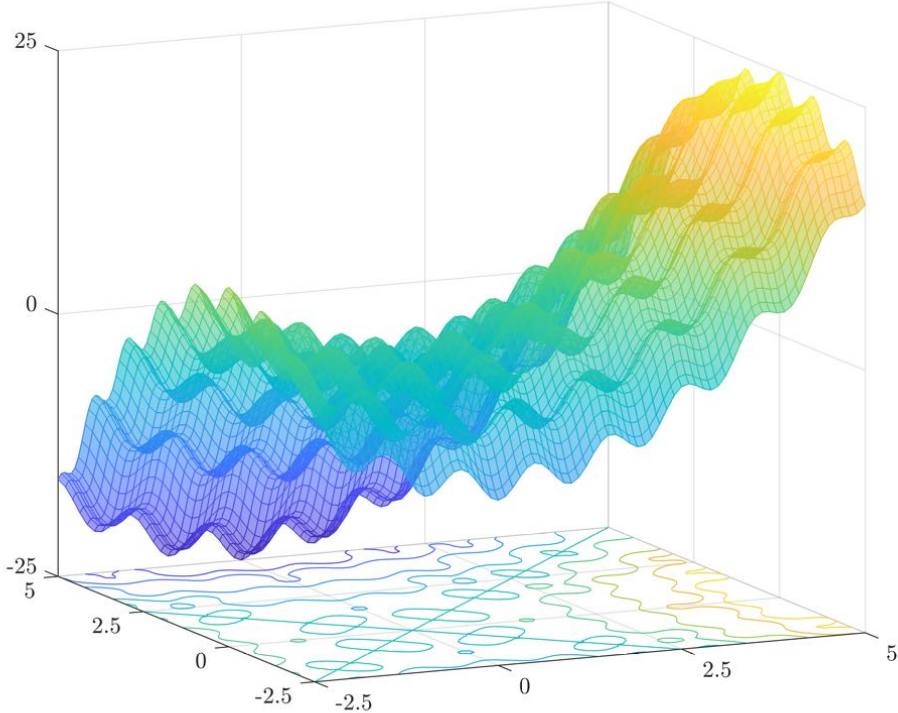}
        \caption{\footnotesize The objective function $\CE(x,y) = R(x)-R(y)$}
    \end{subfigure}~\hspace{\horzspace}~
	\begin{subfigure}[b]{\widthh\textwidth}
        \centering
        \renewcommand{\thesubfigure}{c}
        \includegraphics[trim=85 250 79 250,clip,height=\heightt\textheight]{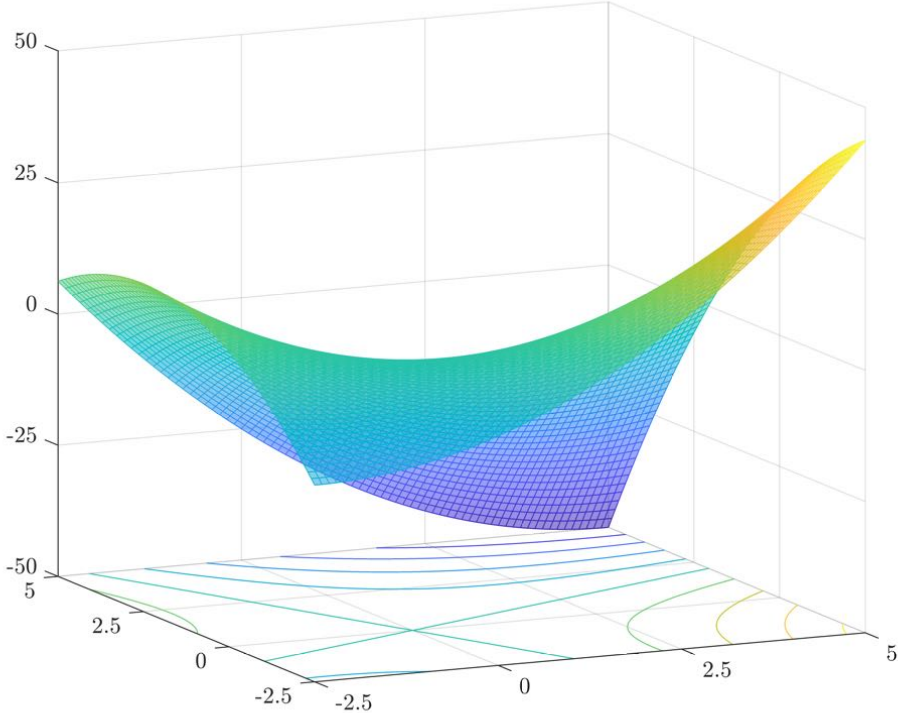}
        \caption{\footnotesize The objective function $\CE(x,y) = x^2-2xy-y^2$}
    \end{subfigure}~\hspace{\horzspace}~
	\begin{subfigure}[b]{\widthh\textwidth}
		\centering
		\renewcommand{\thesubfigure}{d}
        \includegraphics[trim=85 250 79 250,clip,height=\heightt\textheight]{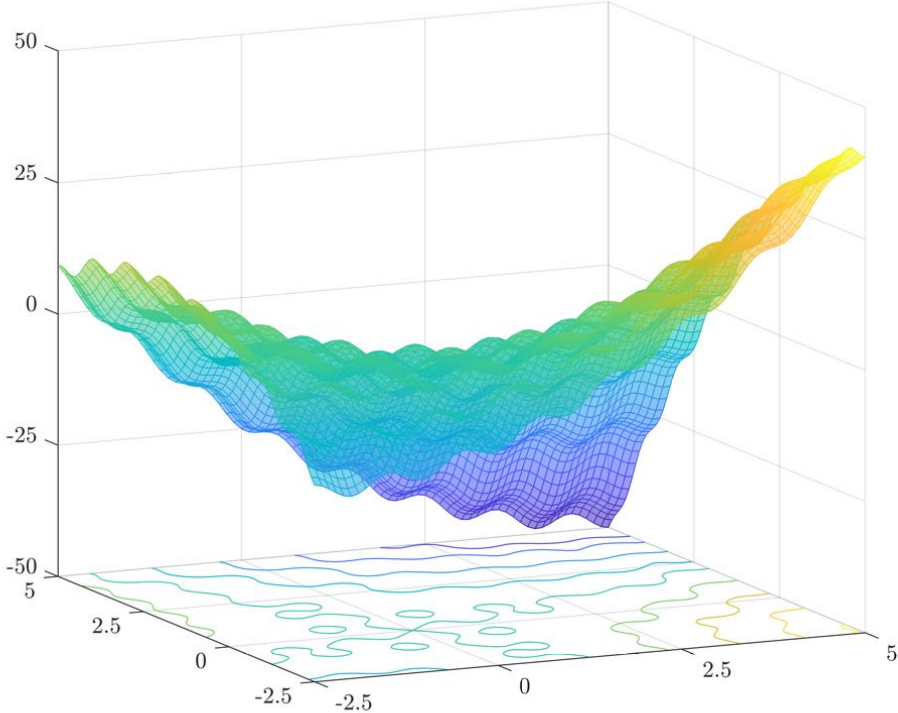}
        \caption{\footnotesize The objective function $\CE(x,y) = R(x)-2xy-R(y)$}
    \end{subfigure}\\\vspace{0.24cm}
    \begin{subfigure}[b]{\widthh\textwidth}
        \centering
        \renewcommand{\thesubfigure}{a0}
        \includegraphics[trim=85 248 92 244,clip,height=\heightt\textheight]{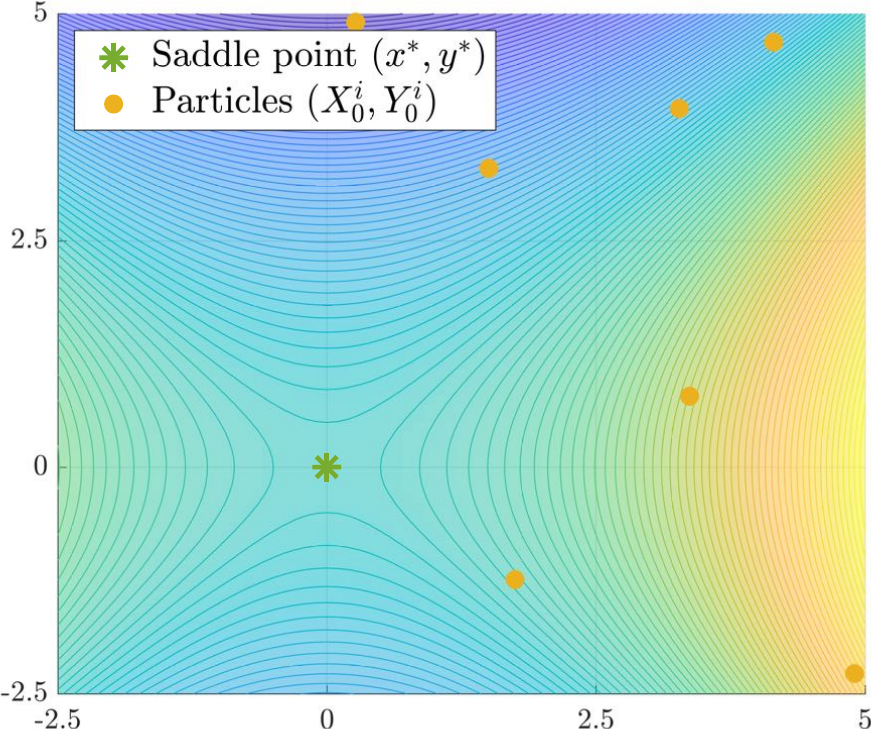}
        \caption{\footnotesize Initial configuration of the particles for \textbf{(a)}}
    \end{subfigure}~\hspace{\horzspace}~
	\begin{subfigure}[b]{\widthh\textwidth}
        \centering
        \renewcommand{\thesubfigure}{b0}
        \includegraphics[trim=85 248 92 244,clip,height=\heightt\textheight]{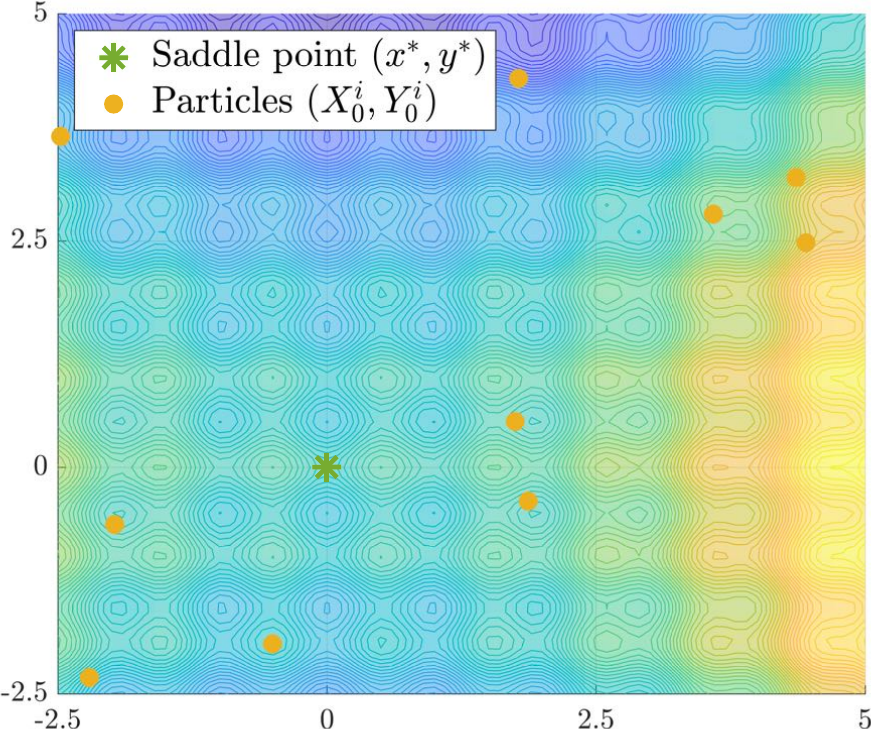}
        \caption{\footnotesize Initial configuration of the particles for \textbf{(b)}}
    \end{subfigure}~\hspace{\horzspace}~
	\begin{subfigure}[b]{\widthh\textwidth}
        \centering
        \renewcommand{\thesubfigure}{c0}
        \includegraphics[trim=85 248 92 244,clip,height=\heightt\textheight]{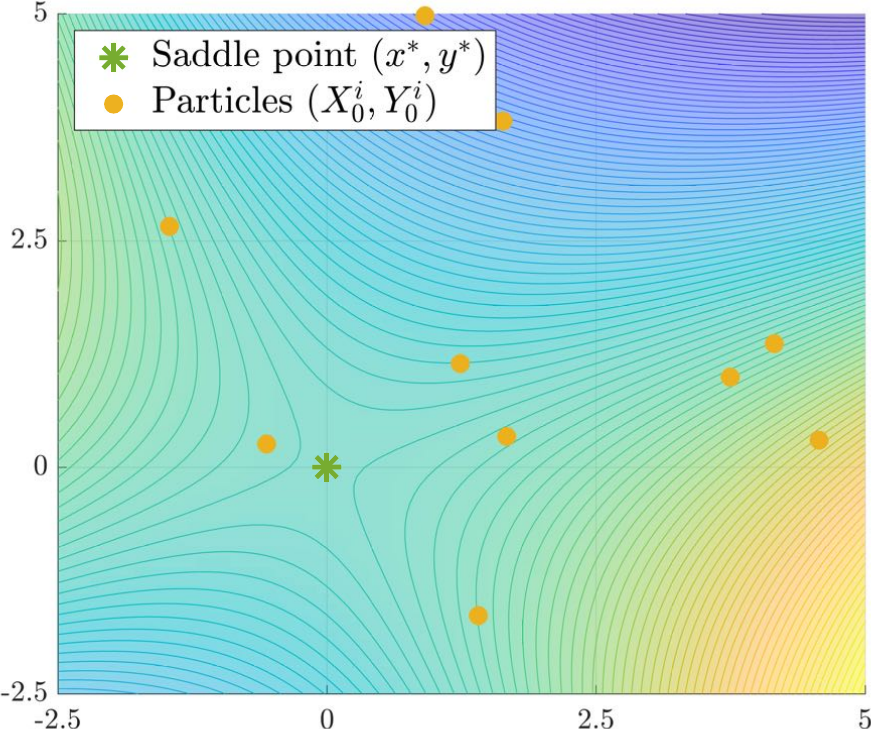}
        \caption{\footnotesize Initial configuration of the particles for \textbf{(c)}}
    \end{subfigure}~\hspace{\horzspace}~
	\begin{subfigure}[b]{\widthh\textwidth}
        \centering
        \renewcommand{\thesubfigure}{d0}
        \includegraphics[trim=85 248 92 244,clip,height=\heightt\textheight]{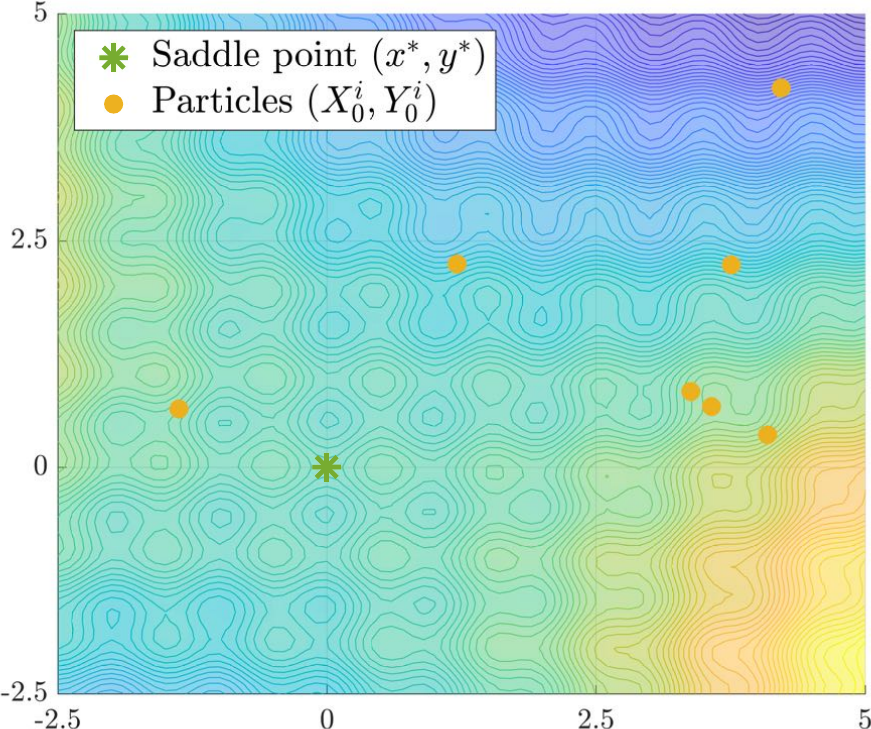}
        \caption{\footnotesize Initial configuration of the particles for \textbf{(d)}}
    \end{subfigure}\\\vspace{0.24cm}
    \begin{subfigure}[b]{\widthh\textwidth}
        \centering
        \renewcommand{\thesubfigure}{a1}
        \includegraphics[trim=85 248 92 244,clip,height=\heightt\textheight]{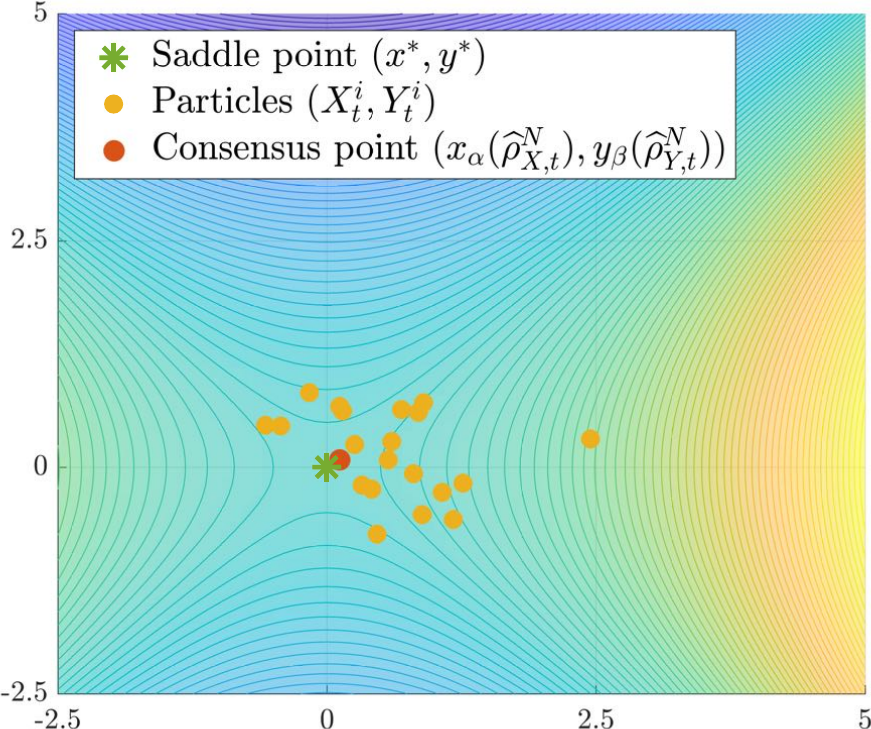}
        \caption{\footnotesize Positions of the particles at $t=2$ for \textbf{(a)}}
    \end{subfigure}~\hspace{\horzspace}~
	\begin{subfigure}[b]{\widthh\textwidth}
        \centering
        \renewcommand{\thesubfigure}{b1}
        \includegraphics[trim=85 248 92 244,clip,height=\heightt\textheight]{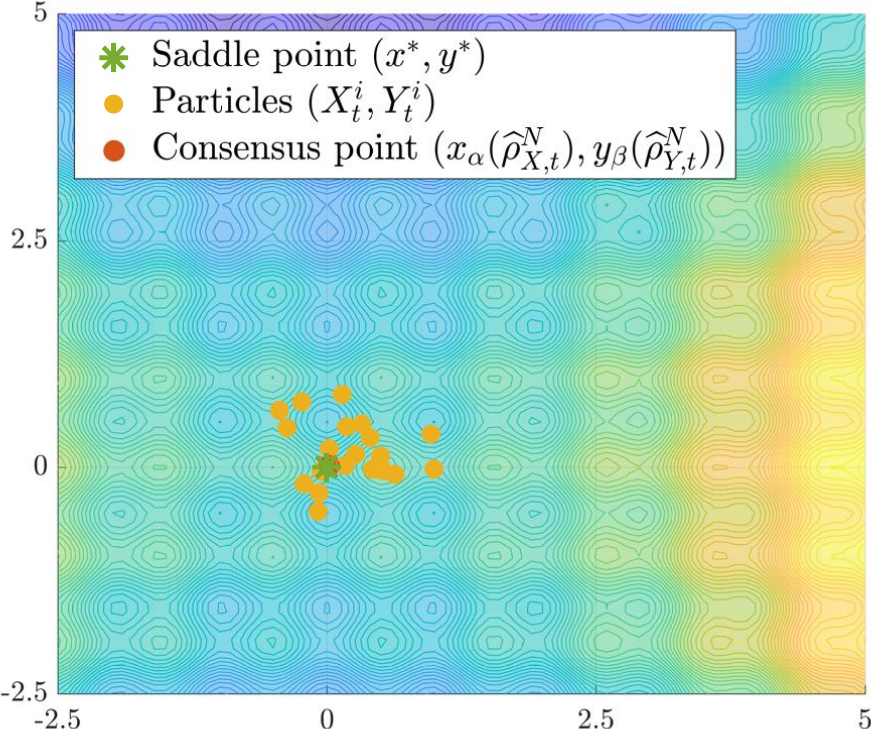}
        \caption{\footnotesize Positions of the particles at $t=2$ for \textbf{(b)}}
    \end{subfigure}~\hspace{\horzspace}~
	\begin{subfigure}[b]{\widthh\textwidth}
        \centering
        \renewcommand{\thesubfigure}{c1}
        \includegraphics[trim=85 248 92 244,clip,height=\heightt\textheight]{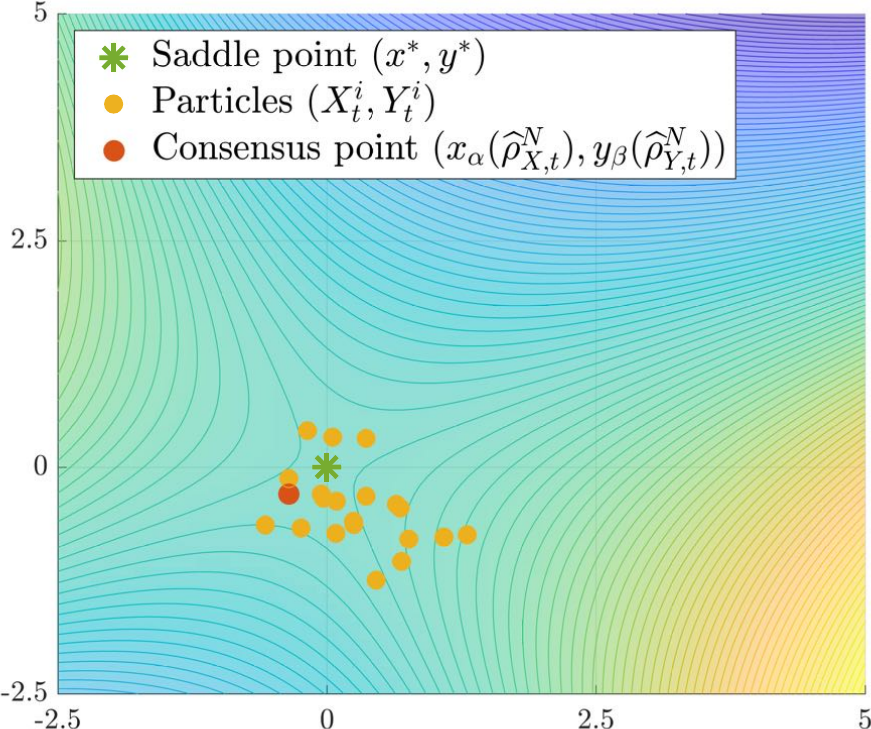}
        \caption{\footnotesize Positions of the particles at $t=2$ for \textbf{(c)}}
    \end{subfigure}~\hspace{\horzspace}~
	\begin{subfigure}[b]{\widthh\textwidth}
        \centering
        \renewcommand{\thesubfigure}{d1}
        \includegraphics[trim=85 248 92 244,clip,height=\heightt\textheight]{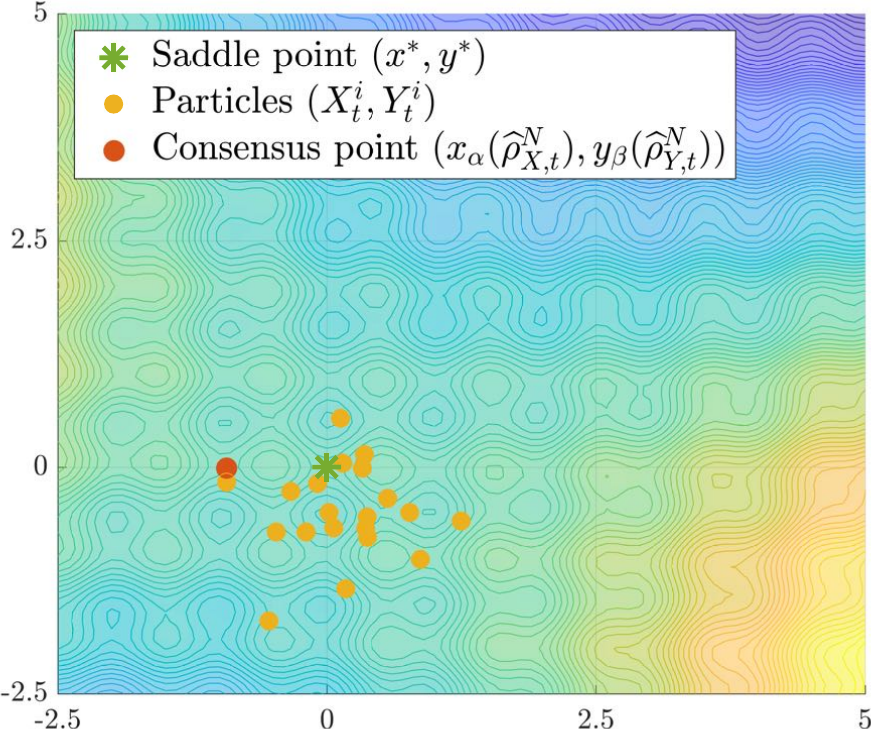}
        \caption{\footnotesize Positions of the particles at $t=2$ for \textbf{(d)}}
    \end{subfigure}\\\vspace{0.24cm}
    \begin{subfigure}[b]{\widthh\textwidth}
        \centering
        \renewcommand{\thesubfigure}{a2}
        \includegraphics[trim=85 248 92 244,clip,height=\heightt\textheight]{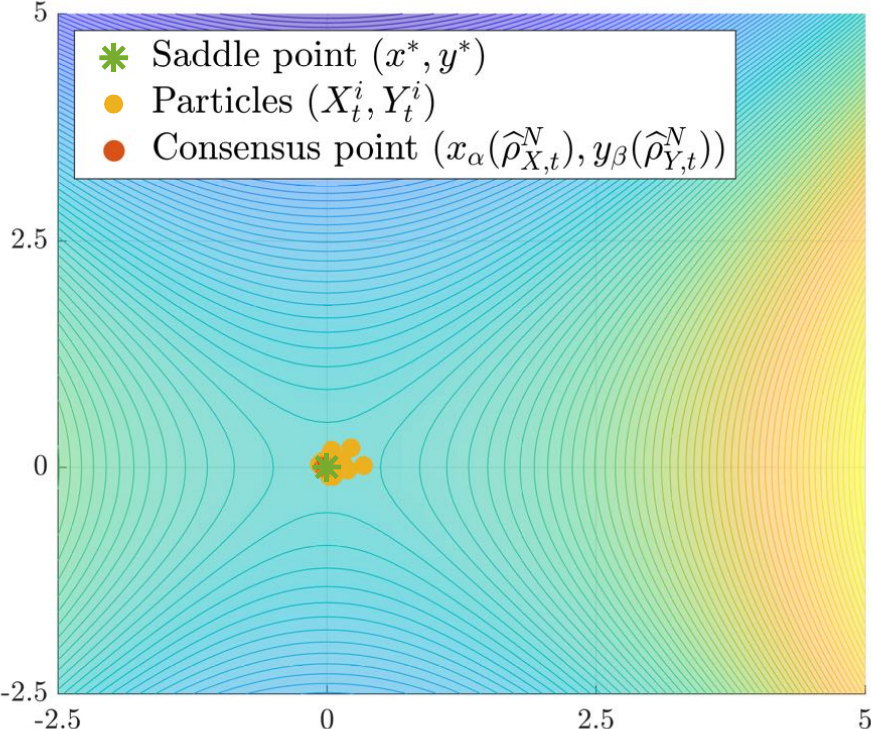}
        \caption{\footnotesize Final configuration of the particles for \textbf{(a)}}
    \end{subfigure}~\hspace{\horzspace}~
	\begin{subfigure}[b]{\widthh\textwidth}
        \centering
        \renewcommand{\thesubfigure}{b2}
        \includegraphics[trim=85 248 92 244,clip,height=\heightt\textheight]{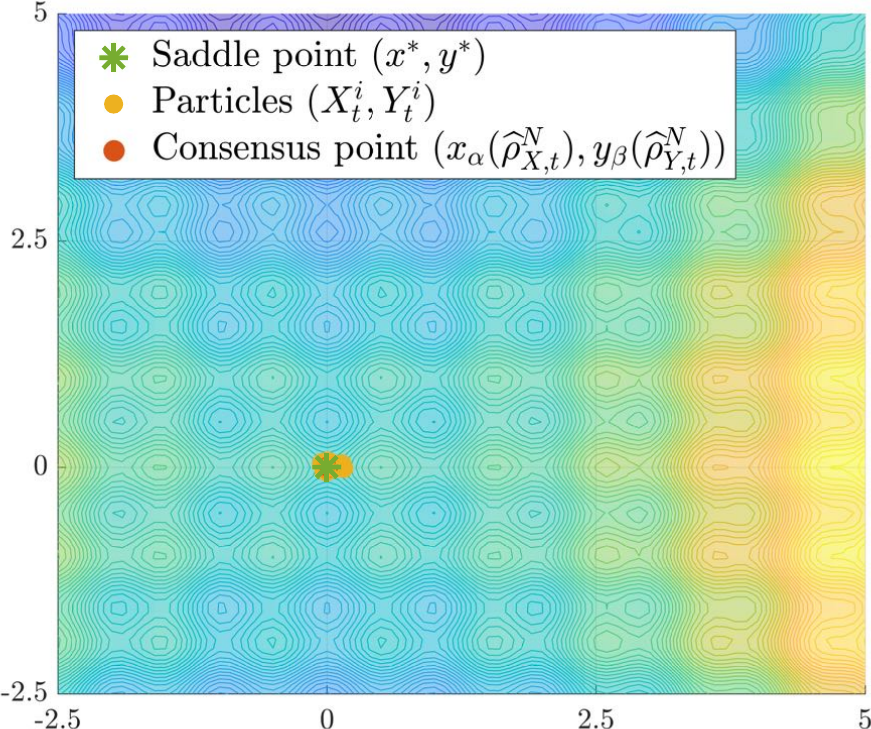}
        \caption{\footnotesize Final configuration of the particles for \textbf{(b)}}
    \end{subfigure}~\hspace{\horzspace}~
	\begin{subfigure}[b]{\widthh\textwidth}
        \centering
        \renewcommand{\thesubfigure}{c2}
        \includegraphics[trim=85 248 92 244,clip,height=\heightt\textheight]{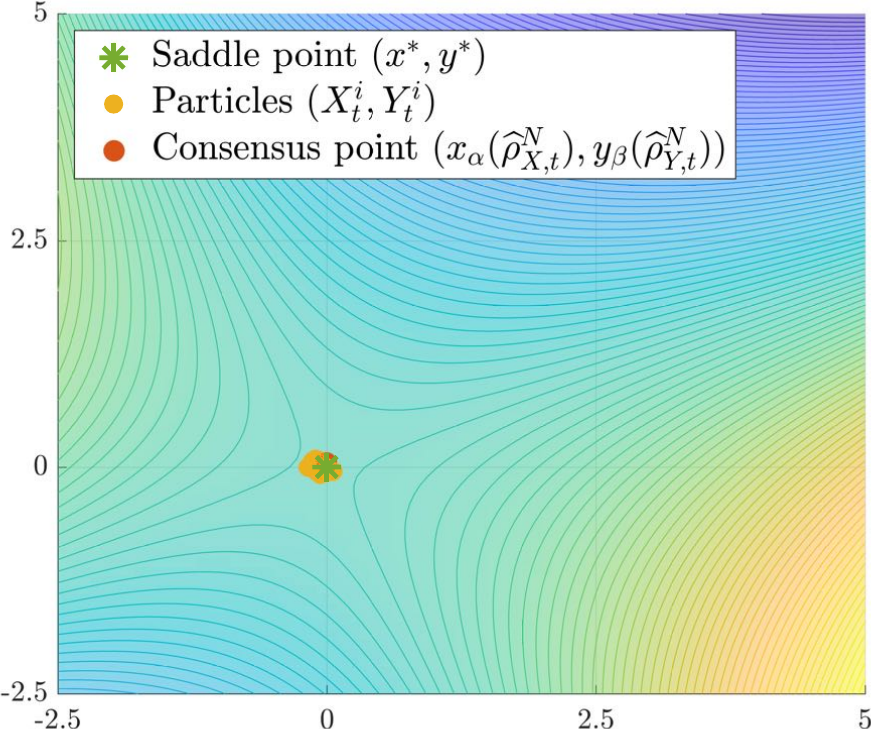}
        \caption{\footnotesize Final configuration of the particles for \textbf{(c)}}
    \end{subfigure}~\hspace{\horzspace}~
	\begin{subfigure}[b]{\widthh\textwidth}
        \centering
        \renewcommand{\thesubfigure}{d2}
        \includegraphics[trim=85 248 92 244,clip,height=\heightt\textheight]{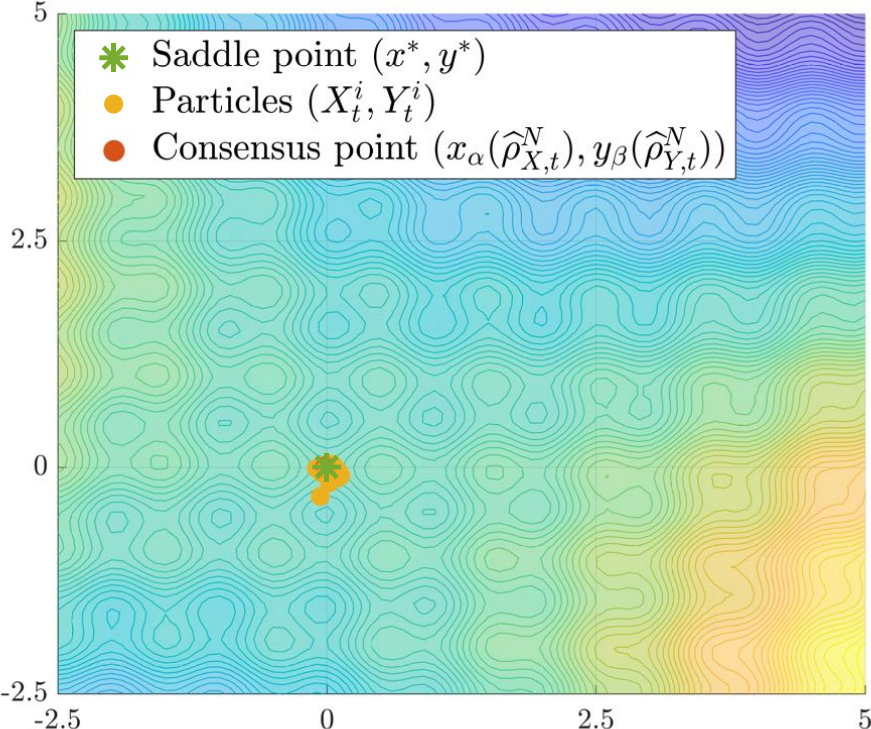}
        \caption{\footnotesize Final configuration of the particles for \textbf{(d)}}
    \end{subfigure}
	\caption{\footnotesize 
	Illustration of the dynamics of CBO-SP when searching the global Nash equilibrium of four different saddle point functions plotted in \textbf{(a)}--\textbf{(d)} in $d=1$, where $R(x) = \sum_{k=1}^d x_k^2 + \frac{5}{2} \big(1-\cos(2\pi x_k)\big)$ is the Rastrigin function. 
	Each column visualizes the positions of the $N=20$ particles when running CBO-SP with parameters~$\alpha=\beta=10^{15}$, $\lambda_1=\lambda_2=1$, $\sigma_1=\sigma_2=\sqrt{0.1}$ and time step size~$\Delta t= 0.1$ at three different points in time ($t=0$, $t=2$ and $t=T=4$).
	The particles are sampled initially from~$\rho_0\sim\CN(2,4)\times\CN(2,4)$.}
	\label{fig:CBOSP_illustrative}
\end{figure}

To visualize the behavior of the CBO-SP algorithm in practice, we depict in Figure~\ref{fig:CBOSP_illustrative} below snapshots of the positions of the particles for four different types of saddle point functions, which are plotted in the first row of the figure.
The experiments include two \textit{nonconvex-nonconcave} examples, which is in general the setting of particular interest in modern applications.
\renewcommand{\floatpagefraction}{.8}
We observe that in all cases (also in case of different initializations) the saddle point is found fast and reliably.

\revised{\subsection{Solving a Quadratic Game with CBO-SP}

To demonstrate the practicability of CBO-SP, we solve a strongly-monotone quadratic game~\cite[Section~5]{loizou2021stochastic} of the form
\begin{align} \label{eq:quadratic_game}
	\min_{x\in\bbR^{d_1}} \max_{y\in\bbR^{d_2}} \frac{1}{n} \sum_{i=1}^n \frac{1}{2}x^TA_ix + x^TB_iy - \frac{1}{2}x^TC_iy
\end{align}
with sample size $n=100$ and for various dimensions $d_1$ and $d_2$.
The matrices $B_i\in\bbR^{d_1\times d_2}$ have random Gaussian entries and the positive definite matrices $A_i\in\bbR^{d_1\times d_1}$ and $C_i\in\bbR^{d_2\times d_2}$ are of the form $A_i = \widetilde{A}_i^T\widetilde{A}_i$ and $C_i = \widetilde{C}_i^T\widetilde{C}_i$ with $\widetilde{A}_i\in\bbR^{d_1\times d_1}$ and $\widetilde{C}_i\in\bbR^{d_2\times d_2}$ having random Gaussian entries.
We employ CBO-SP with parameters~$\alpha=\beta=10^{15}$, $\lambda_1=\lambda_2=1$, $\sigma_1=\sigma_2=2$ using $N\in\{40,80,120,200\}$ particles, and with time horizon $T=100$ and discrete time step size $\Delta t=0.1$.
The particles are sampled initially from~$\rho_0\sim\CN(4,2\Id)\times\CN(4,2\Id)$ (i.e., they are initialized substantially far from the saddle point).
We depict in Table~\ref{table} below the success rates, average $\ell^\infty$-error and average run time of CBO-SP algorithm computed on the basis of $100$ runs.
A run is considered successful if the obtained solution has an accuracy of $10^{-3}$ w.r.t.\@ the $\ell^\infty$-norm.
In brackets, we indicate the average (over the runs) runtime in milliseconds (ms) as well as the average (over the runs) $\ell^\infty$-error.
\begin{table}[ht]
\tiny
	\centering
	\caption{\revised{Success rates, average runtime (in ms) and average $\ell^\infty$-error of the CBO-SP algorithm when solving a quadratic game as specified in \eqref{eq:quadratic_game} for different dimensions $d_1$ and $d_2$, and with different numbers of particles $N$. All results are computed on the basis of $100$ runs of the algorithm.}}
\begin{tabular}[t]{lcccc}
\toprule
& $N=40$ & $N=80$ & $N=120$ & $N=200$\\
\midrule
$d_1=20$, $d_2=8$
     &     $31$\% ($32$ms, $1.5\cdot10^{-2}$)     &     $100$\% ($50$ms, $2.4\cdot10^{-7}$)     &     $100$\% ($182$ms, $6.1\cdot10^{-8}$)     &     $100$\% ($229$ms, $2.9\cdot10^{-8}$)     \\
$d_1=20$, $d_2=20$
     &     $7$\% ($41$ms, $3.0\cdot10^{-2}$)     &     $100$\% ($63$ms, $4.5\cdot10^{-7}$)     &     $100$\% ($339$ms, $3.7\cdot10^{-8}$)     &     $100$\% ($462$ms, $2.4\cdot10^{-8}$)     \\
$d_1=40$, $d_2=8$
     &     $0$\% ($154$ms, $1.1$)     &     $1$\% ($190$ms, $3.6\cdot10^{-2}$)     &     $53$\% ($226$ms, $2.6\cdot10^{-3}$)     &     $100$\% ($418$ms, $4.8\cdot10^{-5}$)     \\
$d_1=40$, $d_2=20$
     &     $0$\% ($162$ms, $1.2$)     &     $0$\% ($285$ms, $4.9\cdot10^{-2}$)     &     $52$\% ($436$ms, $3.8\cdot10^{-3}$)     &     $100$\% ($539$ms, $8.2\cdot10^{-5}$)     \\
$d_1=40$, $d_2=40$
     &     $0$\% ($330$ms, $1.9$)     &     $0$\% ($336$ms, $1.2\cdot10^{-1}$)     &     $25$\% ($421$ms, $5.4\cdot10^{-3}$)     &     $100$\% ($606$ms, $7.9\cdot10^{-5}$)     \\
\bottomrule
\end{tabular}
\label{table}
\end{table}%

\noindent
We observe that with already moderately many particles, the CBO-SP algorithm is capable of consistently finding the desired saddle point for relatively high-dimensional minimax problems.}

Experiments in much higher dimensions and more applied settings coming for instance from economics or arising when training GANs are left to future and more experimental research, which focuses on benchmarking rather than providing rigorous convergence guarantees.

\section{Conclusions} \label{sec:conclusion}

In this paper we propose consensus-based optimization for saddle point problems~(CBO-SP) and analyze its global convergence behavior to global Nash equilibria.
As apparent from the proof, our technique requires the equilibrium to satisfy the saddle point property, i.e., that $\min_x\max_y$ and $\max_y\min_x$ coincide.
We leave to further research the extension of the results to sequential games, where the latter condition does not hold.
This is in particular relevant in, for instance, the training of GANs, which are formulated as non-simultaneous games.

\section*{Acknowledgements}
J.Q.\@~is partially supported by the National Science and Engineering Research Council of Canada (NSERC) and by the PIMS-Europe Fellowship funding.
K.R.\@~acknowledges financial supports from the Technical University of Munich-Institute for Ethics in Artificial Intelligence (IEAI).

The authors moreover sincerely thank the two referees for their careful and insightful comments and suggestions which helped improve the article.

\bibliographystyle{abbrv}
\bibliography{biblio.bib}

\appendix
\revised{
\section{Appendix}

\subsection{Existence and Uniqueness of Solutions to SDEs}

For the sake of a self-consistent presentation, let us recall two results from \cite{Durrett} about the existence and pathwise uniqueness of a strong solution of a SDE of the form
\begin{align} \label{SDE} \tag{$\star$}
	Z_t = Z_0 + \int_0^t b(Z_s) \,ds + \int_0^t\sigma(Z_s) \,dB_s.
\end{align}
These results are used in the proof of Theorem~\ref{theorem:well-posedness_interacting particle dynamcis}.
Note that here we adopted the notation of \cite{Durrett}, i.e., in our setting we have $Z_t=\mathbf{Z}_t$ as well as $b(\mathbf{Z}_t) = -\bm{\lambda}\mathbf{F}(\mathbf{Z}_t)$ and $\sigma(\mathbf{Z}_t) = \bm{\sigma}\mathbf{M}(\mathbf{Z}_t)$.
			
\begin{theorem}[{\cite[Chapter 5, Theorem 3.1]{Durrett}}]\label{thm31}
	Suppose that
	\begin{enumerate}[label=(\roman*),labelsep=10pt,leftmargin=35pt]
		\item\label{thm31_i} for any $n<\infty$ we have $\abs{b_i(z)-b_i(z')}\leq K_n\abs{z-z'}$ and $\abs{\sigma_{ij}(z)-\sigma_{ij}(z')}\leq K_n\abs{z-z'}$ for $\abs{z},\abs{z'} \leq n$,
		\item\label{thm31_ii} there is a constant $A<\infty$ and a function $\varphi(z)\geq 0$ so that if $Z_t$ is a solution of~\eqref{SDE}, then $e^{-At}\varphi(Z_t)$ is a local supermartingale.
	\end{enumerate}
 Then \eqref{SDE} has a strong solution and pathwise uniqueness holds.
\end{theorem}

\begin{theorem}[{\cite[Chapter 5, Theorem 3.2]{Durrett}}] \label{thm32}
	Let $a=\sigma\sigma^T$ and suppose that $\sum_{i=1}^d 2z_i b_i(z) + a_{ii}(z)
		\leq B(1+\abs{z}^2)$.
	Then \ref{thm31_ii} in Theorem \ref{thm31} holds with $A=B$ and $\varphi(z)=1+\abs{z}^2$.
\end{theorem}

\subsection{The Laplace Principle}

\begin{lemma}\label{lmx}
	For any fixed $y\in \bbR^{d_2}$ and any $\delta>0$, define the set
	\begin{align*}
		S_{y,\delta} = \left\{x\in\bbR^{d_1}:\exp\left(-\CE(x,y)\right)>\exp\left(-\min_{x\in\bbR^{d_1}}\CE(x,y)\right)-\delta \right\}.
	\end{align*}
	Let $\{\mu^\alpha\}_{\alpha\geq 1}$ be a family of measures in $\mathcal{P}(\bbR^{d_1})$ and assume that there exists some constant $C_\delta>0$ depending only on $\delta$ such that $\mu^\alpha(S_{y,\delta})\geq C_\delta$ for all $\alpha\geq 1$.
	Then it holds
	\begin{equation*}\label{laplace}
		\lim_{\alpha\to\infty}-\frac{1}{\alpha}\log\left(\int_{\bbR^{d_1}}\exp\left(-\alpha\CE(x,y)\right) d\mu^\alpha(x)\right)
		= \min_{x\in\bbR^{d_1}}\CE(x,y).
	\end{equation*}
\end{lemma}
\begin{proof}
	We first notice that by definition of the set $S_{y,\delta}$ it holds
	\begin{align*}
		\left(\int_{\bbR^{d_1}}\exp\left(-\alpha\CE(x,y)\right) d\mu^\alpha(x)\right)^{\frac{1}{\alpha}}
		&\geq \left(\int_{S_{y,\delta}}  \left(\exp\left(-\textstyle\min_{x\in\bbR^{d_1}}\CE(x,y)\right)-\delta\right)^\alpha d\mu^\alpha(x)\right)^{\frac{1}{\alpha}}\\
		&= \left(\exp\left(-\textstyle\min_{x\in\bbR^{d_1}}\CE(x,y)\right)-\delta\right) \mu^\alpha(S_{y,\delta})^{\frac{1}{\alpha}}\notag\\
		&\geq\left(\exp\left(-\textstyle\min_{x\in\bbR^{d_1}}\CE(x,y)\right)-\delta\right) C_\delta ^{\frac{1}{\alpha}}\\
		&\rightarrow \exp\left(-\textstyle\min_{x\in\bbR^{d_1}}\CE(x,y)\right)-\delta\quad\text{ as }\quad\alpha \to \infty.
	\end{align*}
	Thus, for any $\delta>0$, we have
	\begin{align*}
		\liminf_{\alpha\to\infty}\left(\int_{\bbR^{d_1}}\exp\left(-\alpha\CE(x,y)\right) d\mu^\alpha(x)\right)^{\frac{1}{\alpha}}
		\geq \exp\left(-\textstyle\min_{x\in\bbR^{d_1}}\CE(x,y)\right) - \delta.
	\end{align*}
	On the other hand, clearly
	\begin{align*}
		\limsup_{\alpha\to\infty}\left(\int_{\bbR^{d_1}}\exp\left(-\alpha\CE(x,y)\right) d\mu^\alpha(x)\right)^{\frac{1}{\alpha}}
		\leq  \exp\left(-\textstyle\min_{x\in\bbR^{d_1}}\CE(x,y)\right).
	\end{align*}
	Since $\delta$ was arbitrary, this implies
	\begin{equation*}
		\lim_{\alpha\to\infty}\left(\int_{\bbR^{d_1}}\exp\left(-\alpha\CE(x,y)\right) d\mu^\alpha(x)\right)^{\frac{1}{\alpha}}
		= \exp\left(-\textstyle\min_{x\in\bbR^{d_1}}\CE(x,y)\right),
	\end{equation*}
	giving the result after taking the logarithm on both sides.
\end{proof}
Analogously we obtain the following.
\begin{lemma}\label{lmy}
	For any fixed $x\in \bbR^{d_1}$ and any $\delta>0$, define the set
	\begin{align*}
		S_{x,\delta} = \left\{y\in\bbR^{d_2}:\exp\left(\CE(x,y)\right)>\exp\left(\max_{y\in\bbR^{d_2}}\CE(x,y)\right)-\delta \right\}.
	\end{align*}
	Let $\{\mu^\beta\}_{\beta\geq 1}$ be a family of measures in $\mathcal{P}(\bbR^{d_2})$ and assume that there exists some constant $C_\delta>0$ depending only on $\delta$ such that $\mu^\beta(S_{x,\delta})\geq C_\delta$ for all $\beta\geq 1$.
	Then it holds
	\begin{equation*}\label{laplace}
		\lim_{\beta\to\infty} \frac{1}{\beta}\log\left(\int_{\bbR^{d_2}}\exp\left(\beta\CE(x,y)\right) d\mu^\beta(y)\right)
		= \max_{y\in\bbR^{d_2}}\CE(x,y).
	\end{equation*}
\end{lemma}
}
\end{document}